\documentclass[twoside]{article}
\usepackage{amssymb}
\usepackage[a4paper]{geometry}
\usepackage{a4wide}
\usepackage[latin1]{inputenc} 
\usepackage[T1]{fontenc} 
\usepackage{RR}
\usepackage{hyperref}
\newcommand{\PFeStar}{\mbox{$({\mathcal P}_{\beta}^*)$}}

\newtheorem{theorem}{Theorem}

\newtheorem{lemma}{Lemma}
\newtheorem{proposition}{Proposition}
\newtheorem{remark}{Remark}
\newcommand{\dps}{\displaystyle}
\newcommand{\ssz}[1]{{\scriptsize{#1}}}
\newcommand{\oline}[1]{\dps{\overline{#1}}}
\newcommand{\on}{\mbox{ on }}
\newcommand{\inn}{\mbox{ in }}
\newcommand{\aand}{\mbox{ and }}
\newcommand{\Find}{\mbox{ Find }}
\newcommand{\st}{\mbox{ such that }}
\newcommand{\DF}{\mbox{Darcy-Forchheimer }}
\newcommand{\half}{\scriptsize\frac{1}{2}}
\newcommand{\third}{\scriptsize\frac{1}{3}}
\newcommand{\tthird}{\scriptsize\frac{2}{3}}
\newcommand{\thalf}{\scriptsize\frac{3}{2}}
\newcommand{\fourth}{\scriptsize\frac{1}{4}}

\newcommand{\wklim}{\rightharpoonup}
\newcommand{\emb}{\hookrightarrow}
\newcommand{\lra}{\longrightarrow}
\newcommand{\deff}{:=}
\newcommand{\dfrac}{\displaystyle\frac}
\newcommand{\text}{\mbox}
\newcommand{\bproof}{{\bf Proof: }}
\newcommand{\eproof}{{\hfill $\Box$}\\}

\newcommand{\R}{R}

\newcommand{\dom}{\mathcal O}
\newcommand{\espace}{\dom}
\newcommand{\bord}{\partial\espace}
\newcommand{\dsig}{ds}

\renewcommand{\div}{\mbox{$\mathrm{ div}$}}
\newcommand{\ssdiv}{\!\!\mbox{\scriptsize{ div}}}
\newcommand{\divt}{\mbox{$\mathrm{  div}$}}
\newcommand{\Div}{\,\mbox{$\mathrm{ Div\,}$}}
\newcommand{\grad}{\mbox{$ \nabla$}}
\newcommand{\gradt}{\mbox{$ \nabla$}}
\newcommand{\dsum}[1]{\dps{\sum_{i=1}^{2}{#1}}}
\newcommand{\dint}[1]{\dps{\int_{#1}}}
\newcommand{\dintx}[2]{\dps{\int_{#1}{#2}\;dx}}

\newcommand{\dsup}[2]{\dps{\sup_{#1}{#2}}}

\newcommand{\norm}[2]{\|{#2}\|_{#1}}
\newcommand{\normp}[3]{\|{#2}\|_{#1}^{#3}}
\newcommand{\dpair}[4]{<#3,#4>_{#1,#2}}
\newcommand{\dpairbO}{\dpair{H^{\half}(\partial{\mathcal O})}{H^{-\half}(\partial{\mathcal O})}{\pbo}{\vv\cdot\vn}}
\newcommand{\dpairbg}{\dpair{H^{\half}(\partial\gamma)}{H^{-\half}(\partial\gamma)}{p_{d,\gamma}}{\vvg\cdot\vng}}
\newcommand{\dpairGi}{{{\dpair{H^{\half}(\Gamma_i)}{H^{-\half}(\Gamma_i)}{p_{d,i}}{\vvi\cdot\vni}}}}
\newcommand{\pp}{\mbox{$p$}}
\renewcommand{\pi}{\mbox{$p_i$}}
\newcommand{\pgam}{\mbox{$p_{\gamma}$}}
\newcommand{\pbet}{\mbox{$p_\beta$}}
\newcommand{\pbeti}{\mbox{$p_{\beta,i}$}}
\newcommand{\pbetg}{\mbox{$p_{\beta,\gamma}$}}
\newcommand{\plim}{\mbox{$\tilde{\pp}$}}
\newcommand{\pilim}{\mbox{$\tilde{\pp}_{i}$}}
\newcommand{\pglim}{\mbox{$\tilde{\pp}_{\gamma}$}}

\newcommand{\po}{\mbox{$p_{_{\mathcal O}}$}}
\newcommand{\pbeto}{\mbox{$p_{_{\beta,\mathcal O}}$}}
\newcommand{\pbetjo}{\mbox{$p_{_{\beta_j,\mathcal O}}$}}
\newcommand{\rr}{\mbox{$r$}}

\newcommand{\rgam}{\mbox{$r_\gamma$}}
\newcommand{\qq}{\mbox{$q$}}

\newcommand{\vu}{{\bf{u}}}
\newcommand{\vui}{\vu_i}
\newcommand{\vug}{\vu_\gamma}
\newcommand{\vulim}{\mbox{$\tilde{\vu}$}}
\newcommand{\vuilim}{\mbox{$\tilde{\vu}_i$}}
\newcommand{\vuglim}{\mbox{$\tilde{\vu}_{\gamma}$}}
\newcommand{\vuinclim}{\mbox{$\hat{\vu}_i$}}

\newcommand{\vubet}{\vu_\beta}
\newcommand{\vubeti}{\vu_{\beta,i}}
\newcommand{\vubetg}{\vu_{\beta,\gamma}}
\newcommand{\vuo}{\mbox{$\vu_{_{\mathcal O}}$}}
\newcommand{\vubo}{\mbox{$\vu_{_{\beta,\mathcal O}}$}}
\newcommand{\vubjo}{\mbox{$\vu_{_{\beta_j,\mathcal O}}$}}
\newcommand{\vv}{{\bf{v}}}
\newcommand{\vvi}{\vv_i}
\newcommand{\vvg}{\vv_\gamma}
\newcommand{\vw}{{\bf{w}}}
\newcommand{\vwi}{\vw_i}
\newcommand{\vwg}{\vw_\gamma}

\newcommand{\vn}{{\bf{n}}}
\newcommand{\vni}{{\bf{n}}_{i}}
\newcommand{\vng}{{\bf{n}}_\gamma}

\newcommand{\vx}{{\bf{x}}}
\newcommand{\vy}{{\bf{y}}}
\renewcommand{\gg}{{\bf {g}_{_{\mathcal O}}}}
\newcommand{\ff}{f_{_{\mathcal O}}}

\newcommand{\pbo}{p_{_{\partial,{\mathcal O}}}}

\newcommand{\qo}{q_{_{\mathcal O}}}

\newcommand{\ffo}{f_{_{\mathcal O}}}
\newcommand{\vg}{{\bf g}}

\newcommand{\ali}{\alpha_i}
\newcommand{\beti}{\beta_i}
\newcommand{\algam}{\alpha_\gamma}
\newcommand{\betgam}{\beta_\gamma}
\newcommand{\kgamn}{\kappa}
\newcommand{\alo}{\alpha_{_{\mathcal O}}}
\newcommand{\beto}{\beta_{_{\mathcal O}}}
\newcommand{\ulalo}{ \underline{\alpha}_{_{\mathcal O}}}

\newcommand{\ulalgam}{ \underline{\alpha}_\gamma}
\newcommand{\ulalg}{ \underline{\alpha}_\gamma}
\newcommand{\ulbetgam}{ \underline{\beta}_\gamma}
\newcommand{\ulbetg}{ \underline{\beta}_\gamma}
\newcommand{\ulali}{ \underline{\alpha}_i}

\newcommand{\ulbeti}{ \underline{\beta}_i}

\newcommand{\ulbeta}{ \underline{\beta}}
\newcommand{\ulkap}{\underline{\kgamn}}
\newcommand{\olalo}{ \overline{\alpha}_{_{\mathcal O}}}

\newcommand{\olalgam}{ \overline{\alpha}_\gamma}
\newcommand{\olbetgam}{ \overline{\beta}_\gamma}

\newcommand{\olali}{ \overline{\alpha}_i}
\newcommand{\olalg}{ \overline{\alpha}_\gamma}
\newcommand{\olbeti}{ \overline{\beta}_i}
\newcommand{\olbetg}{ \overline{\beta}_\gamma}
\newcommand{\olalpha}{ \overline{\alpha}}
\newcommand{\olbeta}{ \overline{\beta}}
\newcommand{\olkap}{\overline{\kgamn}}
\newcommand{\tbeta}{\dps{\frac{c_{\ell}^{3}}{2}\beta}}

\newcommand{\tbetgam}{\dps{\frac{c_{\ell}^{3}}{2}\ulbeta_\gamma}}
\newcommand{\hbeta}{C_{\ell}^{3}\beta}

\newcommand{\gbeta}{\frac{c_\ell^3}{2}}

\newcommand{\Cali}{C_L\olali}
\newcommand{\Calg}{C_L\olalg}

\newcommand{\kap}{\kappa}
\newcommand{\xikap}{{\xi}\,{\olkap}}
\newcommand{\xibar}{\bar{\xi}}

\newcommand{\ST}{[\vu_1\cdot\vn -\vu_2\cdot\vn]}

\newcommand{\STbetj}{[\vu_{\beta_j,1}\cdot\vn -\vu_{\beta_j,2}\cdot\vn]}
\newcommand{\FTi}{(\alpha_i +\beta_i |\vu_i|)}
\newcommand{\divugminust}{\hat{\vu}_{\gamma}}
\newcommand{\Fsscript}{\beta}
\newcommand{\D}{\mathcal D}
\newcommand{\Ddef}{ (\D(\oline{\Omega}_1))^{d}\times (\D(\oline{\Omega}_2) )^d\times (\D(\oline{\gamma}))^{d-1}}
\newcommand{\M}{\mathcal M}
\newcommand{\ssM}{_{\mathcal M}}
\newcommand{\Mp}{{\mathcal M}^\prime}
\newcommand{\ssMp}{_{{\mathcal M}^\prime}}

\newcommand{\MFe}{{\mathcal M}_{\Fsscript}}
\newcommand{\ssMFe}{_{{\mathcal M}_{\Fsscript}}}
\newcommand{\MFep}{{\mathcal M}^{\prime}_{\Fsscript}}
\newcommand{\ssMFep}{_{{\mathcal M}^{\prime}_{\Fsscript}}}
\newcommand{\vW}{{\bf{W}}}
\newcommand{\ssW}{_{\bf{W}}}
\newcommand{\vWp}{{\bf W}^\prime}
\newcommand{\ssWp}{_{{\bf W}^\prime}}
\newcommand{\tvW}{\widetilde{\vW}}

\newcommand{\WFe}{{\bf W}_{\!\Fsscript}}
\newcommand{\vWFe}{{\bf W}_{\!\Fsscript}}
\newcommand{\ssWFe}{_{{\bf W}_{\!\Fsscript}}}
\newcommand{\WFep}{\mbox{${\bf W}^{\prime}_{\!\Fsscript}$}}

\newcommand{\ssWFep}{_{{\bf W}^{\prime}_{\!\Fsscript}}}
\newcommand{\tWFe}{\widetilde{\vW}_{\!\Fsscript}}

\newcommand{\vV}{{\bf{V}}}
\newcommand{\ssV}{_{\bf{V}}}
\newcommand{\vVp}{{\bf V}^\prime}
\newcommand{\ssVp}{_{{\bf V}^\prime}}

\newcommand{\VFe}{{\bf V}_{\!\Fsscript}}
\newcommand{\ssVFe}{_{{\bf V}_{\!\Fsscript}}}
\newcommand{\VFep}{\mbox{${\bf V}^{\prime}_{\!\Fsscript}$}}
\newcommand{\ssVFep}{_{{\bf V}^{\prime}_{\!\Fsscript}}}

\newcommand{\WpW}{{\ssWp}_,{\ssW} }
\newcommand{\WFepW}{{\ssWFep}_,{\ssWFe} }
\newcommand{\MpM}{{\ssMp}_,{\ssM} }

\newcommand{\MFepM}{_{{\mathcal M}^{\prime}_{\Fsscript},{\mathcal M}_{\Fsscript}}}
\newcommand{\Omi}{\Omega_{i}}
\newcommand{\spac}{\dom}
\newcommand{\X}{\M_{\beta}({\mathcal O})}

\newcommand{\tX}{\M({\mathcal O})}
\newcommand{\tXp}{\M({\mathcal O})^\prime}
\newcommand{\Y}{\vW_{\!\!\beta}({\mathcal O})}
\newcommand{\tlY}{\widetilde{\vW}_{\!\!\beta}({\mathcal O})}

\newcommand{\tY}{\vW\!({\mathcal O})}
\newcommand{\tltY}{\widetilde{\vW}\!({\mathcal O})}
\newcommand{\tYp}{\vW\!({\mathcal O})^\prime}

\newcommand{\ab}{a_{\beta}}
\newcommand{\ao}{a_{_{\mathcal O}}}
\newcommand{\bo}{b_{_{\mathcal O}}}
\newcommand{\abo}{a_{_{\beta,\mathcal O}}}

\newcommand{\aFe}{a_{{\Fsscript}}}
\newcommand{\AFe}{A_{{\Fsscript}}}
\newcommand{\bFe}{b_{{\Fsscript}}}
\newcommand{\BFe}{B_{{\Fsscript}}}

\newcommand{\Cmap}{C}
\newcommand{\Prob}{\mbox{$({\mathcal P})$}}

\newcommand{\PDar}{\mbox{$({\mathcal P}_{Darcy})$}}
\newcommand{\PForch}{\mbox{$({\mathcal P}_{Forch})$}}

\newcommand{\PFe}{\mbox{$({\mathcal P}_{\Fsscript})$}}
\newcommand{\PFehom}{\mbox{$({\mathcal P}_{\Fsscript}^{0})$}}

\newcommand{\thetab}{\mbox{${\theta_{\beta}}$}}
\newcommand{\theto}{\mbox{${\theta_{_{\mathcal O}}}$}}
\newcommand{\thetbo}{\mbox{${\theta_{_{\beta,\mathcal O}}}$}}
\newcommand{\CLis}{\mbox{${C_{L}^{2}}$}}
\newcommand{\CLgs}{\mbox{${C_{L}^{2}}$}}
\newcommand{\Clic}{\mbox{${C_{\ell}^{3}}$}}
\newcommand{\Clis}{\mbox{${C_{\ell}^{2}}$}}
\newcommand{\Clgc}{\mbox{${C_{\ell}^{3}}$}}
\newcommand{\Clgs}{\mbox{${C_{\ell}^{2}}$}}
\newcommand{\clic}{\mbox{${c_{\ell}^{3}}$}}

\newcommand{\Done}{D_1}
\newcommand{\Dtwo}{D_2}
\newcommand{\Dthree}{D_3}
\newcommand{\Dfour}{D_4}
\newcommand{\Dfive}{D_5}
\newcommand{\Dsix}{D_6}

\newcommand{\Cxi}{C_\xi}
\newcommand{\Cfour}{C_4}
\newcommand{\Cfive}{C_5}
\newcommand{\gORpbord}{\|\gg  \|_{\tYp}}

\RRNo{8443}
\RRdate{December 2013}
\RRauthor{
Peter Knabner\thanks{University of Erlangen-Nuremberg, Department of Mathematics, Cauerstr. 11, D-91058 Erlangen, Germany,
 \hspace*{.8cm} 
e-mail: knabner@am.uni-erlangen.de}%
  \and
Jean E. Roberts\thanks{Inria Paris-Rocquencourt,  B.P. 105, F-78153, Le Chesnay, France,
email: jean.roberts@inria.fr}}%
\authorhead{P. Knabner \& J. E. Roberts}
\RRtitle{Analyse math\'ematique d'un mod\`ele discret de fractures couplant un \'ecoulement de Darcy dans la matrice avec un \'ecoulement de Darcy-Forchheimer dans la fracture}
\RRetitle{Mathematical analysis of a discrete fracture model coupling Darcy flow in the matrix with Darcy-Forchheimer\\ flow in the fracture}
\titlehead{Coupling Darcy and Darcy-Forchheimer flow in a fractured porous medium}
\RRnote{This work was supported by the Franco-German cooporation program PHC Procope}
\RRresume{Nous nous int\'eressons \`a un mod\`ele d'\'ecoulement dans un milieu poreux avec une fracture. Dans ce mod\`ele
l'\'ecoulement dans la fracture est gouvern\'e par la loi de \DF\\ alors que l'\'ecoulement dans la matrice est gouvern\'e par la loi de Darcy.  Nous proposons une formulation variationelle, mixte pour ce mod\`ele et nous d\'emontrons l'existence et l'unicit\'e de la solution.  Pour montrer l'existence nous proposons aussi une formulation analogue pour un mod\`ele
bas\'e sur un \'ecoulement \DF dans tout le domaine.  Nous montrons l'existence et l'unicit\'e de la solution pour ce
deuxi\`eme mod\`ele et montrons que la solution pour le premi\`er mod\`ele est la limite faible de celle du deuxi\`eme mod\`ele quand le
coefficient de Forchheimer dans la matrice tend vers z\'ero.
}
\RRabstract{We consider a model for flow in a porous medium with a fracture in which the flow in the fracture is governed
by the \DF law while that in the surrounding matrix is governed by Darcy's law.  We give an appropriate mixed,
variational formulation and show existence and uniqueness of the solution.  To show existence we give an analogous formulation
for the model in which the \DF law is the governing equation throughout the domain.  We show existence and
uniqueness of the solution and show that the solution for the model with Darcy's law in the matrix is the weak limit
of solutions of the model with the \DF law in the entire domain when the Forchheimer coefficient in the matrix tends
toward zero.
}
\RRmotcle{\'ecoulement en milieu poreux, fractures, \'ecoulement de \DF, solvabilit\'e, r\'egularisation, op\'erateurs monotones}
\RRkeyword{flow in porous media, fractures,  \DF flow, solvability,  regularization, monotone operators}
\RRprojet{Pomdapi}
\RCParis 

\begin{document}
\makeRR   
\section*{Introduction}\label{Intro}
%
%
%
Numerical modeling of fluid flow in a porous medium, even single-phase, incompressible fluid flow, is complicated
because the permeability coefficient characterizing the medium may vary over several orders of magnitude within a region quite small
in comparison to the dimensions of the domain.  This is in particular the case when fractures are present in the medium.   Fractures have
at least one dimension that is very small, much smaller than a reasonable discretization parameter given the size of the domain, but are
much more permeable (or possibly, due to crystalization , much less permeable) than the surrounding medium.   They thus have a very
significant influence on the fluid flow but adapting a standard finite element or finite volume mesh to handle flow in the fractures poses
obvious problems.   Many models have been developed to study fluid flow in porous media with fractures.   Models may employ a
continuum representation of fractures as in the double porosity models derived by homogenization or they may be discrete fracture models.
Among the discrete fracture models  are models of discrete fracture networks in which only the flow in the fractures is considered.  The
more complex discrete fracture models couple flow in the fractures  or in fracture networks with flow in the surrounding medium.  This
later type model is the type considered here.

An alternative to the possibility of using a very fine grid in the fracture and a necessarily much coarser grid away from the fracture is the possibility of treating the fracture as an $(n-1)-$dimensional
hypersurface in the $n-$dimensional porous medium.  This is the idea that was developed in \cite{AlboinJRS[99]} for highly permeable fractures
and in \cite{MartinJR[05]} for fractures that may be highly permeable or nearly impermeable.  Similar models have also been studied in
\cite{FailleEtAl[02],AngotBH[09],MorelesS[10]}.  These articles were all concerned with the case of single-phase, incompressible
flow governed by Darcy's law and the law of mass conservation.  In \cite{FrihRS[08]}  a model was derived in which  Darcy's law was replaced by
the \DF law for the flow in the fracture, while Darcy's law was maintained for flow in the rest of the medium.  The model was
approximated numerically with mixed finite elements and some numerical experiments were carried out.

The use of the linear Darcy law as the constitutive law for fluid flow in porous media, together with the continuity equation, is well established.
For medium-ranged velocities it fits well with experiments \cite[Chapter 5]{Bear[72]} and can be derived rigorously (on simpler periodic media) by homogenization
starting from Stokes's equation \cite{Tartar[80],Allaire[89],Allaire[97]}.  However, for high velocities experiments show deviations which indicate the
need for a nonlinear correction term, \cite{Forchheimer[1901]}, \cite[Chapter 5]{Bear[72]}.  The simplest proposed is a term quadratic in velocity, the Forchheimer correction.
In fractured media, the permeability (or hydraulic conductivity) in the fractures is generally much greater than in the surrounding medium so that the total flow process in the limit is dominated by the fracture flow. 
This indicates that a modeling different from Darcy's model is necessary  
and leads us to investigate models combining Darcy and \DF flow.

In this paper we consider existence and uniqueness of the solution of corresponding stationery problems.  Assumptions on coefficients should
be weak so as not to prevent the use of the results in more complex real life situations.
Therefore we aim at weak
solutions of an appropriate variational formulation, where we prefer a mixed variational formulation, due to the structure of the problems and a
further use of mixed finite element techniques. For a simple d-dimensional domain $\Omega$
and for the linear Darcy flow the
results are well known (c.f. \cite{Brezzi[74]}) and rely on the coercivity of the operator $A$ coming from Darcy's equation on the kernel of the divergence
operator $B$ coming from the continuity equation and the functional setting in $H(\div,\Omega)$ for the flux and $L^2(\Omega)$ for the pressure.
For the nonlinear \DF flow the functional setting has to be changed to $W^3(\div,\Omega)$ (see Appendix \ref{Wpdivspaces}) for the flux so that $A$ will remain
(strictly) monotone and to  $L^{\thalf}(\Omega)$ for the  pressure.  This makes it possible to extend the reasoning for the linear case to the
homogeneous \DF problem and via regularization, using the Browder-Minty theorem for maximal monotone operators, also to prove unique
existence in the inhomogeneous case.  This work is carried out in the thesis \cite{Summ-th[01]};  see also \cite{Knabner-Summ, Fabrie[89], Amirat[91]}
for related results.  Here we extend this reasoning to the situation of two subdomains of the matrix separated by a fracture with various choices
of the constitutive laws in domains and fractures.  One would expect that the \DF law is more accurate than Darcy's law (and this will be partially
made rigorous); therefore, (and for technical reasons) we start with a model having the \DF law
 throughout the domain (though with strongly variable coefficients) and extend the aforementioned
reasoning for existence and uniqueness to this case, (Section \ref{ExUniqDFFrac+DFMatrix}).
By its derivation, Darcy's law should be a limit case of the \DF law.  This is made precise in
Section 4 by showing that the solution of the Darcy model is a weak limit of solutions of the \DF
model with the Forchheimer coefficient (multiplying the nonlinear term) going to 0. This was shown earlier in \cite{Amirat[91]} under slightly different assumptions,
but we include it here for completeness.  This opens up the possibility of treating
various combinations of the constitutive laws.  As rapid transport is more likely to take place  in the fractures, we explicitly treat the case of Darcy's law
in the matrix and  the \DF law in fractures.  By using the full \DF model as a regularization and deriving corresponding a priori bounds we can
show the existence of a solution as a weak limit of the regularizing full models (Section \ref{ExUniqDFFrac+DMatrix}).  Uniqueness again follows as in all the other cases
from the monotone structure of the problem (see Appendix \ref{ApUniqueness}).  Technical difficulties stem from the different functional settings for the linear case
and the nonlinear case.  It may be envisaged to extend this basic procedure in various directions.  An obvious extension is to
the case of a finite number of fractures and subdomains, as long as the fractures do not intersect, which is quite
restrictive.  But also a general case where d-dimensional subdomains are separated by (d-1)-dimensional fractures,
which are separated by (d-2)-dimensional fractures, etc. may be attacked with this approach.  Another extension could be the investigation of other nonlinear correction terms to Darcy's law: cf. \cite{BMikelicW[09]}.

The outline of this article is as follows:  in  Section \ref{FormOfProbs} the model problem with Darcy flow in the
matrix and \DF flow in the fracture as well as the problem with \DF\\ flow in the matrix and in the fracture will
be given.  In Section \ref{ExUniqDFFrac+DFMatrix} the existence and uniqueness of the solution to the problem with \DF flow
in the matrix and in the fracture will be shown.  Section \ref{SimpDom} is concerned with showing that in a simple domain (one without a fracture)
that the solution of the Darcy problem is obtained as the limit of the \DF problem
when the Forchheimer coefficient tends to zero.  Then Section \ref{ExUniqDFFrac+DMatrix} takes up the problem for extending the result of Section \ref{SimpDom} to the case of a domain with a fracture in which it is shown
that the problem with \DF flow in the fracture but with Darcy flow in the matrix is obtained
as the limit of the problem with \DF flow everywhere as the Forchheimer coefficient in the
matrix tends to zero.
%
%
\section{Formulation of the problems}\label{FormOfProbs}   
%
%
\subsection{Formulation with \DF flow in the fracture and Darcy flow in the matrix}\label{DFFrac+DMatrix}  
%
Let $\Omega$ be a bounded domain in $\R^d$ with boundary $\Gamma$, and let $\gamma\subset\Omega$ be a
$(d-1)$-dimensional surface that separates $\Omega$ into two subdomains:
$\Omega\subset \R^d,\quad\Omega=\Omega_1\cup\gamma\cup\Omega_2, \quad \gamma=(\overline{\Omega}_1\cap \overline{\Omega}_2)\cap\Omega,
\quad\Gamma=\partial\Omega, \aand \quad \Gamma_i=\Gamma\cap\partial\Omega_i$.
We suppose for simplicity that $\gamma$ is a subset of a hyperplane; i. e. that $\gamma$ is flat. Taking the stratification of natural porous media into account this seems to be a feasible assumption covering a variety of situations. The extension to the case that $\gamma$ is  a smooth surface should not pose any major problems but would be considerably more complex as the curvature tensor would enter into the definitions of the tangential gradient and the tangential divergence.
We consider the following problem, which was derived in \cite{FrihRS[06], FrihRS[08]}:
\begin{equation}\label{eqnsubd}\begin{array}{rlll}
	\alpha_i \vu_i +\grad \pi &=&0 &\quad\inn\Omega_i\\
	\div\,\vu_i&=&q_i&\quad\inn\Omega_i\\
	p_i&=&p_{d,i} &\quad\on\Gamma_i\\
\end{array}\end{equation}
together with 	
\begin{equation}\label{eqnfrac}\begin{array}{rlll}
	(\algam+\betgam |\vug|)\vug +\gradt\pgam&=&0 &\quad\on\gamma\\
	\divt\,\vug&=&q_\gamma  + \ST &\quad\on\gamma\\
	\pgam&=&p_{d,\gamma} &\quad\on\partial\gamma
\end{array}\end{equation}
and the interface condition
\begin{equation}\label{eqncoup}\begin{array}{rlll}\hspace*{2cm}	
	\pi &=& \pgam +(-1)^{i+1}\kap(\xi\vu_i\cdot\vn+\xibar \vu_{i+1}\cdot\vn),\qquad\on\gamma,
	\quad i=1,2,
\end{array}\end{equation}
where $\vn$ is the unit normal vector on $\gamma$, directed outward from $\Omega_1, \;\kap$ is a coefficient function on $\gamma$ related directly to the fracture width and inversely to the normal component of the permeability of the physical fracture,
the parameter $\xi$ is a constant greater than $1/2$ and $\xibar=1-\xi$, and for convenience of notation the index $i$ of the subdomains
 is considered to be an element of $Z_2$ (so that if $i=2$, then $i+1=1)$.
The tensor coefficients $\alpha_i, i=1,2, \aand \alpha_\gamma$ are related to the inverse of the permeability tensors on $\Omega_i, i=1,2, \aand \gamma,$ respectively, and the  coefficient $\beta_\gamma$ is the Forchheimer coefficient on $\gamma$, assumed to be scalar.
We assume that  the functions
$\;\ali\,:\,\Omega_i \longrightarrow R^{d,d},
	\qquad \algam\,:\,\gamma \longrightarrow R^{d-1,d-1},$
	are all symmetric and uniformly positive definite:
\begin{equation}\label{hypcoef}\begin{array}{rccllll}
	\ulali  |\vx|^2 &\leq& \vx\cdot\ali(\vy) \vx &\leq&\olali |\vx|^2
		&\quad \forall \vy\in\Omega_i,&\;\vx\in R^d\\
	\ulalgam |\vx|^2 &\leq& \vx\cdot\algam(\vy) \vx &\leq&{\olalgam}  |\vx|^2
		&\quad \forall \vy\in\gamma,&\;\vx\in R^{d-1}\\
   \text{and}\qquad\quad    \ulbetgam &\leq& \betgam(\vy) &\leq&{\olbetgam} 
	&\quad \forall \vy\in\gamma,
\end{array}\end{equation}	
where $\quad \ulali,\,  \ulalgam, \, \ulbetgam >0$,  and that the real valued coefficient
function $\kap \,:\,\gamma \longrightarrow R$ is bounded above and below by
positive constants:
\begin{equation}\label{hypkap}
	0<\ulkap\leq\kap(\vy)\leq\olkap\quad\forall \vy\in \gamma.
\end{equation}
Note that only minimal assumptions concerning $\ali,\,i=1,2,  \algam$ and $\betgam$ reflecting the structure of the problem
 are required and no further regularity, allowing for general heterogeneous media.  However this means that the standard functional setting
 of the linear case has to be modified  and thus also the regularity requirements concerning the source and boundary terms.

We make the following assumptions concerning the data functions $q$ and $p_d$
corresponding respectively to an external source term and to Dirichlet boundary data:
\begin{equation}\label{hypdata}\begin{array}{c}
	q=(q_1,q_2,q_{\gamma})\in L^{{3}}(\Omega_1)\times L^{{3}}(\Omega_2) \times L^{3}(\gamma)\\
	p_d=(p_{d,1},p_{d,2},p_{d,\gamma})\in
				( W^{\third,\thalf}(\Gamma_1)\cap H^{\half}(\Gamma_1))
		\times (W^{\third,\thalf}(\Gamma_2)\cap H^{\half}(\Gamma_2))
		\times H^{\half}(\partial\gamma),
\end{array}\end{equation}
where we have used the standard notation for the Lebesgue spaces $L^p,\,p\in\R,\,p\geq 1,$ and for the Sobolev spaces $W^{k,p},\; k,p\in\R,\,p\geq 1;$ see \cite{Adams[75]}.
Following standard practice we often write $H^k$ for the Sobolev space $W^{k,2},\,k\in\R.$
We have required more regularity of the data functions than necessary for a weak
formulation of problem (\ref{eqnsubd}), (\ref{eqnfrac}), (\ref{eqncoup}) in order to use
the same data functions for problem (\ref{eqnsubd}), (\ref{eqnfrac}), (\ref{eqncoup})
and for problem (\ref{eqnsubdFe}), (\ref{eqnfracFe}), (\ref{eqncoupFe}) given below.

To give a weak mixed formulation of  problem  (\ref{eqnsubd}), (\ref{eqnfrac}), (\ref{eqncoup}), we introduce several spaces of functions:
\[\begin{array}{c}
	\M=\{p=(p_1,p_2,\pgam) \,:\, \pi    \in L^{2}(\Omega_i), i=1,2, \aand \pgam\in L^{3/2}(\gamma)\} \\
		\|p\|\ssM  =
			\displaystyle{\sum_{i=1}^{2}\|p_i\|_{0,2,\Omega_i}
				\,+\, \|\pgam\|_{0,\frac{3}{2},\gamma} }.
\end{array}\]
The space $\M$ being a product of reflexive Banach spaces is clearly a reflexive Banach space with the dual space
\[\begin{array}{c}
	\Mp=\{f=(f_1,f_2,f_\gamma) \,:\, f_i   \in L^2(\Omega_i), i=1,2, \aand  f_\gamma\in L^{3}(\gamma)\} \\
		\|f\|\ssMp  =
			\displaystyle{\sum_{i=1}^{2}\|f_i\|_{0,2,\Omega_i}\,+\, \|f_\gamma\|_{0,3,\gamma} }.
\end{array}\]
We also define
\[\begin{array}{c}
	\vV=\{\vv=(\vv_1,\vv_2,\vvg) \,:\, \vv_i    \in (L^2(\Omega_i))^d, i=1,2,
		\aand \vvg\in (L^{3}(\gamma))^{d-1}  \}  \\
	\|\vv\|\ssV  =
			\displaystyle{\sum_{i=1}^{2}\|\vv_i\|_{0,2,\Omega_i}\,+\, \|\vvg\|_{0,3,\gamma}   }
\end{array}\]
and its dual space
\[\begin{array}{c}
	 \vVp=\dps{\{\vg=(\vg_1,\vg_2,\vg_\gamma) \,:\, \vg_i    \in (L^{2}(\Omega_i))^d, i=1,2,
		\aand \vg_\gamma\in (L^{\frac{3}{2}}(\gamma))^{d-1}}  \}     \\
	\|\vg\|\ssVp  =
			\displaystyle{\sum_{i=1}^{2}\|\vg_i\|_{0,{2},\Omega_i}\,+\, \|\vg_\gamma\|_{0,\frac{3}{2},\gamma}    }.
\end{array}\]
\begin{remark}\label{lonenorm}  For $f=(f_1,f_2,f_\gamma)\in\Mp$, respectively
$\vg=(\vg_1,\vg_2,\vg_\gamma)\in\vVp$,
we have used the $\ell^1$ norm on $R^3$ to
give the norm of $f$, respectively $\vg$, in terms of its three components $f_1,f_2 \aand f_\gamma$,
respectively $\vg_1,\vg_2 \aand \vg_\gamma$, whereas the
actual norm for the dual space would have used the $\ell^\infty$ or maximum norm.  However
these norms are equivalent since $R^3$ is of finite dimension and we have found it more convenient to use the $\ell^1$ norms here.  	
\end{remark}	

For the domains $\Omega_1, \Omega_2$ in $\R^d$ and $\gamma$ in $\R^{d-1}$, respectively, we need minimal regularity to make some of the expressions used below
well defined.  In particular, we need exterior normal vector fields on the boundaries.  To assume that the domains are Lipschitzian will be sufficient, and this
will be done henceforth.
We will need in addition the space $\vW$ defined by
\begin{equation}\label{defW}\begin{array}{c}
	\vW=\{\vu=(\vu_1,\vu_2,\vug)\in \vV\,:\,\Div\vu \deff (\div\vu_1,\div\vu_2, \divt \vug -\ST)\in \Mp\\
			\aand \vu_i\cdot\vn\in L^2(\gamma),  i=1,2  \}\\
	\|\vu\|\ssW =  \|\vu\|\ssV \,+\,\|\Div\vu\|\ssMp\,+\,
		\displaystyle{\sum_{i=1}^{2}\,
			\,\|\vu_i\cdot\vn\|_{0,2,\gamma})	}.			
\end{array}\end{equation}
One can show that $\vW$ is also a reflexive Banach space and that
\begin{equation} \label{DefD}  \D=\Ddef  \end{equation}
is dense in $\vW$ (see e.g. \cite[Lemma 3.13]{Summ-th[01]}), where by $\D(\oline{\dom})$ is meant $\{\psi_{|\dom}:\psi\in C^\infty (\R^{n})\}$,
for $\dom$ a bounded domain in $\R^{n}$.
We also have that
for $\vv\in\vW, \vvi\in H(\div,\Omi),i=1,2, \aand \vvg\in H(\divt,\gamma) \;$(since $L^3(\gamma)\subset
L^2(\gamma)$) so that $\vvi\cdot\vni \in H^{-\half}(\partial\Omi)$ and $\vvg\cdot\vng
\in H^{-\half}(\partial\gamma)$, where $\vni,\;i=1,2,$ and $\vng$ are the exterior normal vectors
on $\partial\Omi,\; i=1,2,$ and on $\partial\gamma$, respectively.

{{Define the  forms}} $a\,:\, \vW\times\vW\longrightarrow \R$ and $b\,:\, \vW\times\M\longrightarrow \R$ by
\[\begin{array}{lllll}
	a(\vu,\vv) = \displaystyle{\sum_{i=1}^{2}\,\int_{\Omega_i}\ali \vu_i\cdot\vv_i\,dx
		+\int_{\gamma}(\algam+\betgam|\vug|)\vug\cdot\vvg \,\dsig
		+\sum_{i=1}^{2}\int_{\gamma}\kap(\xi\vu_i\cdot\vn+\xibar\vu_{i+1}\cdot\vn)
		\vv_i\cdot\vn\,\dsig},\\[.3cm]
		b(\vu,r) = \displaystyle{\sum_{i=1}^{2}\,\int_{\Omega_i}\div\vu_i r_i\,dx
		+\int_{\gamma}(\divt\vug
		- \ST\,\,) \rgam\,\dsig} \;=\;<\Div\vu,r>\MpM .
\end{array}\]
{Note that the form $a$ is continuous and linear in its second variable while the form $b$ is clearly continuous and bilinear.}
Define the continuous, linear forms $\vg\in \vWp$ and $f\in \Mp$  by
\begin{equation}\begin{array}{lllll}\label{Def-f-g}
&\vg\,:\,\vW\longrightarrow \R \\
&\vg(\vv)=-\displaystyle{\sum_{i=1}^{2}\dpairGi\, -\dpairbg}\\[.4cm]
 \mbox{ and }& f\,:\,\M\longrightarrow \R \\
	& f(r)=\displaystyle{\sum_{i=1}^{2}\,\int_{\Omega_i} q_{i} r_i \,dx
	 	+\int_{\gamma}q_\gamma\rgam \, \dsig}.
\end{array}\end{equation} 
The weak mixed formulation of (\ref{eqnsubd}), (\ref{eqnfrac}) and (\ref{eqncoup}) is
	\[ \Prob\quad
			\begin{array}{l}
			\Find \vu\in\vW \aand p\in\M \st\\[.1cm]
			\begin{array}{llllcllll}
				a(\vu,\vv)- b(\vv,p) &=& \vg(\vv) &\quad \forall \vv\in\vW\\[.1cm]
				b(\vu,r) 			  &=& f(r) 	&\quad \forall r\in\M  .\\
			\end{array}
			\end{array}\]
Define also,  for the moment only formally (see Lemma \ref{AFeContSMono}),
\[\begin{array}{lllll}
	A\,:\,\vW\longrightarrow \vWp  &\quad\aand \quad&  B\,:\,\vW\longrightarrow \Mp\\
	<A(\vu),\vv>\WpW = a(\vu,\vv)\quad\forall \vv\in\vW   &&
			<B(\vu),r>\MpM = b(\vu,r)	\quad\forall r\in\M
\end{array}\]
and note that $B:\vW\lra\Mp$ is simply $\Div:\vW\lra\Mp$ so that for $\tvW,$ the kernel of $B,$
\[\begin{array}{c}
\tvW=\{\vv=(\vv_1,\vv_2,\vvg)\in \vW\,:\,\Div\vv = B(\vv)=0\},
\end{array}\]
we have that					
\[\|\vv\|\ssW = \|\vv\|\ssV  \;+\; \dsum{\|\vv_i \cdot\vn}\|_{0,2,\gamma} \qquad\forall \vv\in\tvW.\]
%
\subsection{Formulation with \DF flow in the matrix and in the fractures}\label{DFFrac+DFMatrix}  
%
With $\Omega, \gamma, \Omega_i, \Gamma_i,\vni,\;i=1,2,\vng$ and $\vn$ as well as $\ali, \algam,\betgam, \kap, q_i,
q_{\gamma}, p_{d,i}, p_{d,\gamma},\xi,$ and $ \xibar$ as in the preceding paragraph, and with
$\beti:\Omi\lra \R$ a function satisfying
\begin{equation}\label{hypbeti}
	\ulbeti  \;\leq\; \beti(\vy) \;\leq\;\olbeti 
		\;\quad \forall\, \vy\in\Omega_i,
\end{equation}
where $\ulbeti ,\olbeti >0$,
we now consider the following problem:
\begin{equation}\label{eqnsubdFe}\begin{array}{rlll}
	\FTi \vu_i +\grad \pi &=&0 &\quad\inn\Omega_i\\
	\div\,\vu_i&=&q_i&\quad\inn\Omega_i\\
	p_i&=&p_{d,i} &\quad\on\Gamma_i
\end{array}\end{equation}
together with
\begin{equation}\label{eqnfracFe}\begin{array}{rlll}
	(\algam+\betgam |\vug|)\vug +\gradt\pgam&=&0 &\quad\on\gamma\\
	\divt\,\vug&=&q_\gamma  + \ST &\quad\on\gamma\\
	\pgam&=&p_{d,\gamma} &\quad\on\partial\gamma
\end{array}\end{equation}
and the interface conditions
\begin{equation}\label{eqncoupFe}\begin{array}{rlll}	
	\pi &=& \pgam +(-1)^{i+1}\kap(\xi\vu_i\cdot\vn+\xibar \vu_{i+1}\cdot\vn),\quad i=1,2.
\end{array}\end{equation}
Due to the Forchheimer regularization in the matrix equations, the spaces in the earlier definitions need to
be replaced by spaces appropriate for the functional setting of the Forchheimer equations, i.e. $L^2(\Omi)$
by $L^{\thalf}(\Omi),$ and consequently $L^2(\Omi)$ by $L^3(\Omi)$ for the dual spaces, thus obtaining
a $\beta$-version of the earlier spaces, i. e. $\MFe$ instead of $\M$, etc.  For the sake of clarity we state
explicitly:
\[\begin{array}{c}
	\MFe=\{p=(p_1,p_2,\pgam) \,:\, \pi    \in L^{3/2}(\Omega_i), i=1,2, \aand \pgam\in L^{3/2}(\gamma)\} \\
		\|p\|\ssMFe  =
			\displaystyle{\sum_{i=1}^{2}\|p_i\|_{0,\thalf,\Omega_i}
				\,+\, \|\pgam\|_{0,\thalf,\gamma} }.
\end{array}\]
The space $\MFe$  is clearly a reflexive Banach space with dual space
\[\begin{array}{c}
	\MFep=\{f=(f_1,f_2,f_\gamma) \,:\, f_i   \in L^3(\Omega_i), i=1,2, \aand  f_\gamma\in L^{3}(\gamma)\} \\
		\|f\|\ssMFep  =
			\displaystyle{\sum_{i=1}^{2}\|f_i\|_{0,3,\Omega_i}\,+\, \|f_\gamma\|_{0,3,\gamma} }.
\end{array}\]

We also define
\[\begin{array}{c}
	\VFe=\{\vv=(\vv_1,\vv_2,\vvg) \,:\, \vv_i    \in (L^3(\Omega_i))^d, i=1,2,
		\aand \vvg\in (L^{3}(\gamma))^{d-1}  \}  \\
	\|\vv\|\ssVFe  =
			\displaystyle{\sum_{i=1}^{2}\|\vv_i\|_{0,3,\Omega_i}\,+\, \|\vvg\|_{0,3,\gamma}   },
\end{array}\]
which is similarly a reflexive Banach space, with its dual space
\[\begin{array}{c}
	 \VFep=\dps{\{\vg=(\vg_1,\vg_2,\vg_\gamma) \,:\, \vg_i    \in (L^{3/2}(\Omega_i))^d, i=1,2,
		\aand \vg_\gamma\in (L^{3/2}(\gamma))^{d-1}}  \}     \\
	\|\vg\|\ssVFep  =
			\displaystyle{\sum_{i=1}^{2}\|\vg_i\|_{0,\thalf,\Omega_i}\,+\, \|\vg_\gamma\|_{0,\thalf,\gamma}.    }
\end{array}\]
Again we have used the equivalent $\ell^1$ norm instead of the $\ell^\infty$ norm to construct the
product space norm for $\MFep\aand \WFep$.
   We also need the space $\WFe$ defined by
\begin{equation}\label{defWFe}\begin{array}{c}
	\WFe=\{\vu=(\vu_1,\vu_2,\vug)\in \VFe\,:\,\Div\vu = (\div\vu_1,\div\vu_2, \divt \vug -\ST)\in \MFep\\
			\aand \vu_i\cdot\vn\in L^2(\gamma),  i=1,2 \}\\
	\|\vu\|\ssWFe =  \|\vu\|\ssVFe \,+\,\|\Div\vu\|\ssMFep\,+\,
		\displaystyle{\sum_{i=1}^{2}\,
			\,\|\vu_i\cdot\vn\|_{0,2,\gamma}.	}			
\end{array}\end{equation}
One can show that $\WFe$ is  a reflexive Banach space, that $\D$, given by (\ref{DefD}), is dense in $\WFe$, that for $\vv\in\WFe,$
for $i=1,2, \vvi\in H(\div,\Omi) \aand \vvg\in H(\divt,\gamma) $
so that $\vvi\cdot\vni \in H^{-\half}(\partial\Omi)$ and
$\vvg \cdot\vng\in H^{-\half}(\partial\gamma)$.  Further $\vvi\in W^{3}(\div,\Omi)$ (see Appendix \ref{Wpdivspaces} ).
{{Define the  forms}} $\aFe\,:\, \WFe\times\WFe\longrightarrow \R$ and $\bFe\,:\, \WFe\times\MFe\longrightarrow \R$ by
\[\begin{array}{lllll}
	\aFe(\vu,\vv) = \displaystyle{\sum_{i=1}^{2}\,\int_{\Omega_i}(\ali+\beti|\vu_i|)\vu_i\cdot\vv_i\,dx
		+\int_{\gamma}(\algam+\betgam|\vug|)\vug\cdot\vvg \,\dsig}\\[.3cm]
		\hspace*{3cm}
		+\dps{\sum_{i=1}^{2}\int_{\gamma}\kap(\xi\vu_i\cdot\vn+\xibar\vu_{i+1}\cdot\vn)\vv_i\cdot\vn\,\dsig},\\[.3cm]
		\bFe(\vu,r) = \displaystyle{\sum_{i=1}^{2}\,\int_{\Omega_i}\div\vu_i r_i\,dx
		+\int_{\gamma}(\divt\vug
		- \ST\,\,) \rgam\,\dsig} \;=\;<\Div\vu,r>\MFepM.
\end{array}\]
Note that the form $\aFe$ is continunous and linear in its second variable while  $\bFe$ is continuous and bilinear.
Define the linear forms
\[\begin{array}[t]{lllll}
	\vg\,:\, \WFe  \longrightarrow \R			&\;\aand\;&         f\,:\,\MFe\longrightarrow \R\\
\end{array}\]
as in (\ref{Def-f-g}) but with $\vg\in \WFep$ and $f\in\MFep$ which is valid with the regularity assumptions in (\ref{hypdata}).
The mixed weak formulation of (\ref{eqnsubdFe}), (\ref{eqnfracFe}) and (\ref{eqncoupFe}) is given by
	\[ \PFe\quad
			\begin{array}{l}
			\Find \vu\in\WFe \aand p\in\MFe \st\\[.1cm]
			\begin{array}{llllcllll}
				\aFe(\vu,\vv)- \bFe(\vv,p) &=& \vg(\vv) &\quad \forall \vv\in\WFe\\[.1cm]
				 \bFe(\vu,r) 			&=& f(r) 	&\quad \forall r\in\MFe.\\
			\end{array}
			\end{array}\]

Define again
\[\begin{array}{lllll}\hspace*{.6cm}
	\AFe\,:\,\WFe\longrightarrow \WFep  &\aand \quad&  \BFe\,:\,\WFe \longrightarrow \MFep\\
	<\AFe(\vu),\vv>{\ssVp,\ssV} = \aFe(\vu,\vv)\quad\forall \vv\in\WFe   &&
	\hspace*{-.6cm}
			<\BFe(\vu),r>{\ssMFep,\ssMFe} = \bFe(\vu,r)	\quad\forall r\in\MFe
\end{array}\]
for an equivalent operator equation and
\[\begin{array}{c}
	\tWFe=\{\vu=(\vu_1,\vu_2,\vug)\in \WFe\,:\,\Div\,\vu\deff\BFe(\vu)=0\},
\end{array}\]
and note that
\begin{equation}\label{normWFetilde}\|\vu\|\ssWFe = \|\vu\|\ssVFe  \;+\; \dsum{\|\vu_i \cdot\vn}\|_{0,2,\gamma} \qquad\forall \vu\in\tWFe.\end{equation}
\begin{remark}\label{parambeta}
Note that none of the spaces $\WFe, \VFe, \tWFe,$ or $\MFe$ and neither of the operators $\bFe$ nor $\BFe$
depends on the coefficient  $\beta$.  The index $\beta$ is used simply to indicate that these are the
spaces and operators used to define the problem $\PFe$.
\end{remark}
{To obtain some of the estimates that we will derive in the following sections we shall make
use of the following technical lemma given in \cite[lemmas 1.1 and 1.4]{Knabner-Summ}.}

\begin{lemma} For $\vx$ and $\vy$ in $\R^n$, we have the following inequalities:
\begin{eqnarray}
	\label{ineqa}	|\, |\vx|\vx - |\vy|\vy\, |	&\leq& (|\vx|+|\vy|)\,|\vx-\vy|,\\
	\label{ineqb}  	\half\, |\vx-\vy|^{3}		&\leq&   (|\vx|\vx - |\vy|\vy)\cdot(\vx-\vy) ,  \\
	\label{ineqc}	\left| \, |\vx|^{-\half}\vx - |\vy|^{-\half}\vy\, \right| &\leq& \sqrt{2}\, |\vx-\vy|^{\half},\\
	\label{ineqd}	\dfrac{|\vx-\vy|^{2}}{\sqrt{|\vx|}+\sqrt{|\vy|}}     &\leq&
						\left(\dfrac{\vx}{\sqrt{|\vx|}}-\dfrac{\vy}{\sqrt{|\vy|}}\right)\cdot(\vx-\vy).
\end{eqnarray}
\end{lemma}
In (\ref{ineqc}) and hereafter $x\mapsto |x|^{-\half}x$ on $\R^n$ means the continuation of this function
on $R^{n}\setminus\{0\}$ to $\R^n$ obtained by defining  $|0|^{-\half}0\deff 0$, which by (\ref{ineqc})
is indeed H\"older continuous with exponent $\half$.

Here we introduce some notation that we will use throughout the remainder of the article:  for any
positive integer $n$ and any bounded domain $\spac$ in $\R^n$, we know that
$L^3(\spac)\emb L^2(\spac)$ and that the inclusion map is continuous so that there is a constant
$C_{L,\spac}$ depending on $n$ and the measure of the space such that if $\phi\in L^3(\spac)$
then $\|\phi\|_{L^2}\leq  C_{L,\spac} \|\phi\|_{L^3}$.
Here we shall assume that $C_L$ is a constant with $C_{L,\spac}\leq C_L$ for all of the
spaces $\spac$ that we deal with. (There are only a finite number for each problem.)
  Also we know that if $s \aand t$ are such that
$1\leq s\leq t\leq \infty$
then the $\ell^s\aand\ell^{t}$ norms on $\R^n$ are equivalent (since all norms on finite
dimensional spaces are equivalent), and we shall assume that there are positive real numbers
$C_\ell\aand c_\ell$ such that if $\vx\in\R^n$ then
$c_{\ell}\|\vx\|_{\ell^t}\leq\|\vx\|_{\ell^s}\leq C_\ell\|\vx\|_{\ell^t}$ for all dimensions $n$ and all
norms $\ell^s\aand\ell^t$ with
$1\leq s\leq t\leq \infty$, that we encounter in the problems that follow.  (Again there will only
be a finite number.)
%
%
\section{Existence and uniqueness of the solution of the problem $\PFe$ \DF flow in
the fracture and in the subdomains}\label{ExUniqDFFrac+DFMatrix}			
%
%
To show the existence and uniqueness of the solution of $\PFe$, following 
the argument of \cite[Section 1]{Knabner-Summ} we show that the operator $\AFe\,:\,\WFe\longrightarrow \WFep$ is
continuous and monotone and is uniformly monotone on $\tWFe$ to obtain
a solution to the homogeneous problem with $f=0$.
(That  $\BFe\,:\,\WFe \longrightarrow \MFep$ satisfies the inf-sup condition follows just as in the linear case, cf. \cite{MartinJR[05]}, however
for completeness a demonstration is given in Appendix \ref{ApInfSup}).
Then taking any solution to the second equation of $\PFe$ (whose existence is guaranteed by the inf-sup condition) an auxiliary  homogeneous problem is constructed whose solution can be used to produce the solution of $\PFe$.

%
\begin{lemma}\label{AFeContSMono}
	The operator $\AFe\,:\,\WFe\longrightarrow \WFep$ is continuous and strictly monotone
	and is furthermore uniformly monotone on $\tWFe$.	
\end{lemma}
\bproof{{To see that $\forall \vu\in\WFe, \;\AFe(\vu)\in\WFep\;$ i. e. that $\forall \vu\in\WFe, \;\AFe(\vu)$
			is bounded,}} suppose  that $\vu\in \WFe$.  Then, using the equivalence of norms
			in a finite dimensional space  and H\"older's inequality, we have, for each $\vv$ in $\WFe$,
			\begin{equation}\label{estNonLinTerm} \begin{array}{llllll}
				\dint{\Omega_i }{\beti|\vu_i |\vu_i\cdot\vv_i }\,dx
				&\leq& \olbeti \Clic \dint{\Omega_i}{|\vu_i|^2_3 |\vv_i|_3} dx
				\;\leq\; \olbeti \Clic \displaystyle{
		     		\left(\int_{\Omega_i}(|\vu_i|^2_3)^{\frac{3}{2}}\dsig\right)^{\frac{2}{3}}
				\left(\int_{\Omega_i} |\vv_i|^3_3\dsig\right)^{\frac{1}{3}}}\\[.3cm]
				&\leq& \olbeti \Clic \|\vu_i\|^2_{0,3,\Omega_i} \|\vv_i\|_{0,3,\Omega_i}.
			\end{array}\end{equation}
We also have									
			\begin{equation}\label{estLinTerm} \begin{array}{llllll}
				\dint{\Omega_i }{\ali\vu_i\cdot\vv_i }\,dx
				 \leq  \olali \|\vu_i\|_{0,2,\Omega_i} \|\vv_i\|_{0,2,\Omi}
				 \leq \olali \CLis\|\vu_i\|_{0,3,\Omega_i} \|\vv_i\|_{0,3,\Omega_i}.
			\end{array}\end{equation}
		As we have similar inequalities for the norms on $\gamma$ and as $\xi\geq\xibar$
		we conclude that for each $\vv\in\WFe,$

	\begin{equation}\label{boundAFe} \begin{array}{lll}
		\|\AFe(\vu)\|\ssWFep
		&=&\dps{\sup_{\vv\in\WFe\atop \vv\neq 0}}
				\frac{|<\AFe(\vu),\vv>\WFepW| }{\|\vv\|\ssWFe}\\
		&\leq&\dps{\sup_{\vv\in\WFe\atop \vv\neq 0}}
			\left\{\dsum{\Big( \olali \CLis\|\vu_i\|_{0,3,\Omega_i}
					{ \|\vv_i\|_{0,3,\Omega_i}}
			+\olbeti \Clic \|\vu_i\|^2_{0,3,\Omega_i}
					{ \|\vv_i\|_{0,3,\Omega_i}} \Big)}\right.\\
			&&\qquad\;+\;\dps{\olalgam \CLgs\|\vug\|_{0,3,\gamma}
					{ \|\vvg\|_{0,3,\gamma}}
			+\olbetgam \Clgc \|\vug\|^2_{0,3,\gamma}
					{ \|\vvg\|_{0,3,\gamma}}}\\
			&&\left.\left.\qquad\dps{\;+{\xikap   }}\;
				\left(\dsum{  \|\vu_i\cdot\vn\|_{0,2,\gamma} }\right)
				\left({\dsum{ \|\vv_i\cdot\vn\|_{0,2,\gamma}  }}\right)\right\}\right/{\|\vv\|\ssWFe}\\
	       &\leq&\dsum{\Big( \olali \CLis\|\vu_i\|_{0,3,\Omega_i}
			+\olbeti \Clic \|\vu_i\|^2_{0,3,\Omega_i}    \Big)}\\
	       		&&\qquad\;+\;\dps{\olalgam \CLgs\|\vug\|_{0,3,\gamma}
			+\olbetgam \Clgc \|\vug\|^2_{0,3,\gamma}}
				\dps{\;+\;\xikap   }\;
			\dsum{  \|\vu_i\cdot\vn\|_{0,2,\gamma} }\\
		&\leq& \olalpha C^2_L  \|\vu\|\ssVFe+\olbeta C^3_\ell  \|\vu\|\ssVFe^2
			\;+\;\xikap\;    \dsum{  \|\vu_i\cdot\vn\|_{0,2,\gamma} }\\[.4cm]
		 &\leq& \big(\max\{\olalpha C^2_L  \,,\,  \xikap\}
		 	\;+\; \olbeta C^3_\ell  \|\vu\|\ssVFe  \big) \|\vu\|\ssWFe,

	\end{array}\end{equation}
	where $\olalpha $ is the max\{$\olalpha_1, \olalpha_2, \olalgam$\}, and similarly for $\olbeta$.


{{ To see that $\AFe\,:\,\WFe \longrightarrow\WFep$ is continuous}}
 suppose that $\vu$ and $\vw$ are elements of $\WFe$.  Using H\"older's inequality and
 then inequality (\ref{ineqa}) along with the equivalence of norms in finite dimensional spaces, we see that, for any $\vv\in\WFe$
 and for $i=1,2,$
 \[ \begin{array}{lllll}
	 \dint{\Omi}{\beti( |\vui|\vui-|\vwi|\vwi)\cdot\vvi}\,dx
	\leq  \olbeti \|\, |\vui|\vui-|\vwi|\vwi \, \|_{0,\thalf,\Omega_i} \|\vvi\|_{0,3,\Omega_i} \\[.5cm]
		\leq  \olbeti \Clic\|\, (\|\vui\|_{0,3,\Omega_i}+\|\vwi\|_{0,3,\Omega_i})
			\|\vui-\vwi \|_{0,3,\Omega_i}    \|\vvi\|_{0,3,\Omega_i}.
 \end{array}\]
Then using the analogous inequality for the nonlinear term on $\gamma$ we have
 	\[ \begin{array}{lllll}
	\|\AFe(\vu) \;-\; \AFe(\vw)\|\ssWFep
	&=& \dps{\sup_{\ssz{\begin{array}{c}\vv\in\WFe\\ \vv\neq 0\end{array}}}}
	\dfrac{< \AFe(\vu) - \AFe(\vw)\, , \,\vv >\WFepW}{\|\vv\|\ssWFe}\\[.6cm]
	&\leq& \biggl(\left(  \olalpha  C_L^2  \|\vu-\vw\|\ssVFe
					\;+\; \olbeta  C_\ell^3 (\|\vu\|\ssVFe +\|\vw\|\ssVFe) \|\vu-\vw\|\ssVFe   \right)
					\|\vv\|\ssVFe      \\[.5cm]	
	&&	\hspace{1.3cm}        \;+\; \xikap\;\dsum{\|(\vu_i-\vw_i)\cdot\vn\|_{0,2,\gamma}}
							\dsum{\|\vv_i\cdot\vn\|_{0,2,\gamma}}    \biggr) \biggr/{\|\vv\|\ssWFe}\\[.5cm]
	&\leq& \Bigl( \dps{\max\{\olalpha C^2_L  \,,\,  \xikap\}}
	 		\;+\;  \olbeta  C_\ell^3 (\|\vu\|\ssVFe +\|\vw\|\ssVFe)\Bigr) \|\vu-\vw\|\ssWFe.
	 \end{array}\]

	
{ {To see that $\AFe\,:\,\WFe \longrightarrow\WFep$ is strictly monotone}}
suppose again that $\vu$ and $\vw$ are elements of $\WFe$.  Then using inequality (\ref{ineqb}), for $i=1,2,$
\begin{equation}\label{lowestNonLinTerm}\begin{array}{lllll}
	\dint{\Omega_i}{\beti(|\vu_i|\vu_i-|\vw_i |\vw_i)\cdot(\vu_i -\vw_i)}\,dx
		\geq \dfrac{ \ulbeti \clic}{2}\|\vu_i-\vw_i\|_{0,3,\Omega_i}^3.
\end{array}\end{equation}

We also note that if $x,y\in \R$ then $\xi(x^2+y^2) +2\xibar xy \geq \min\{1,2\xi-1\}(x^2+y^2)$.
It follows that 						
 $\forall \vu, \vw \in \WFe$
	\begin{equation}\label{monotoneAFe}\begin{array}{rll}
		< \AFe(\vu) - \AFe(\vw)\, , \,\vu-\vw >\WFepW \\
		 \geq& \dps{\frac{\ulbeta \clic}{2} \,\|\vu-\vw\|^{3}\ssVFe  }
		+\ulkap\min\{1,2\xi-1\}\dsum{  \|(\vu_i-\vw_i)\cdot\vn\|_{0,2,\gamma}^2}	\\
		\geq& {\mathcal C}(\ulbeta, \ulkap, \xi)\left( \|\vu-\vw\|^{3}\ssVFe
		\;+\; \dsum{  \|(\vu_i-\vw_i)\cdot\vn\|_{0,2,\gamma}^2 }\right)\\
		\geq& 0,
	\end{array}\end{equation}
	where $\ulbeta=\min\{\ulbeta_1,\ulbeta_2,\ulbetg\}$,  and where we have equality only if $\vu=\vw$.


{{To see that $\AFe$ is uniformly monotone on $\tWFe$}} it suffices to note that if $\vu \aand \vw $ belong to $ \tWFe$ then
	\[\begin{array}{lll}
		< \AFe(\vu) - \AFe(\vw)\, , \,\vu-\vw >\WFepW&\geq&
			{\mathcal G}(\|\vu-\vw\|\ssWFe)\|\vu-\vw\|\ssWFe
	\end{array}\]
with
\[\begin{array}{lllllllllllll}
 	{\mathcal G}(\|\vu\|\ssWFe)  &:=&
	 {\mathcal C}(\ulbeta, \ulkap, \xi)\;\dps{\frac{ \|\vu\|^{3}\ssVFe
		\;+\; \dsum{  \|\vui\cdot\vn\|_{0,2,\gamma}^2 }}{\|\vu\|\ssWFe}}\\
	&=& {\mathcal C}(\ulbeta, \ulkap, \xi)\;
	\dps{\frac{ \|\vu\|^{3}\ssVFe \;+\; \dsum{  \|\vui\cdot\vn\|_{0,2,\gamma}^2 }}
		{ \|\vu\|\ssVFe \;+\; \dsum{  \|\vui\cdot\vn\|_{0,2,\gamma} }  }}&\longrightarrow&\infty\quad
			\mbox{ as }\|\vu\|\ssWFe&\longrightarrow&\infty,
\end{array}\]
where we have used (\ref{normWFetilde}). \eproof

\begin{lemma}\label{LeminfsupFe}
The linear form $\bFe\,:\,\WFe\times\MFe\lra\R$ satisfies the following inf-sup condition: there is a positive constant $\thetab$ such that $\forall r\in\MFe$
\begin{equation}\label{infsupFe}
\dps{\thetab \|r\|\ssMFe \leq\sup_{\vv\in{\ssz{\WFe}}}\frac{\bFe(\vv,r)}{\|\vv\|\ssWFe}}.
\end{equation}
\end{lemma}
\bproof
See Appendix \ref{ApInfSup}.
\eproof
\begin{proposition}
	The homogeneous problem
	\[ \PFehom\quad
		\begin{array}{l}
		\Find \vu^{0}_{\Fsscript}\in\WFe \aand p^{0}_{\Fsscript}\in \MFe \st\\
		\begin{array}{llllcllll}
			\aFe(\vu^{0}_{\Fsscript},\vv)&-& \bFe(\vv,p^{0}_{\Fsscript}) &=& \vg(\vv) &\qquad \forall \vv\in\WFe\\
			 \bFe(\vu^{0}_{\Fsscript},r) 	&&		&=& 0 	&\qquad \forall r\in\MFe
		\end{array}
		\end{array}\]
	has a unique solution.\label{prop:prop1}
\end{proposition}
\bproof  That there is a unique solution in $\tWFe$ to $\aFe (\vu^0_{\Fsscript},\vv)=\vg(\vv), \forall\vv\in \tWFe$, i. e.  to $\AFe(\vu^0_{\Fsscript})=\vg,$
now follows from the Browder-Minty theorem, \cite[Theorem 26.A]{Zeidler[90]}.  That there is a
unique $ p^{0}_{\Fsscript}\in \MFe$ such that  $(\vu^{0}_{\Fsscript},p^{0}_{\Fsscript})$ is the
unique solution of $\PFehom$ then follows as in the linear case as the operator $\BFe$
is still linear. 
\eproof\\[.2cm]
To handle a source term in the continuity equation we start from any solution to this equation and construct an auxiliary homogeneous problem whose solution is then combined with the solution to the (nonhomogeneous) continuity equation to produce the desired solution to the full problem. 
\begin{theorem}\label{unisolPFe}
The problem $\PFe$ admits a unique solution $(\vubet,\pbet)\in\WFe\times\MFe$.
\end{theorem}
\bproof
Since, according to Lemma \ref{LeminfsupFe}, $\bFe$ satisfies the inf-sup condition,
 the subproblem of $\PFe$
\[ \begin{array}{l}
            \Find \vu \in\WFe \st\\[.1cm]
             \;b_\beta(\vu,r) = f(r) \quad \forall r \in \M_\beta
 \end{array}\]
has a (non-unique) solution.  Let  $\vu^* \in \vW_\beta$ denote one such.
We consider the auxiliary problem
\[ \PFeStar \quad
            \begin{array}{l}
            \Find \vulim\in\WFe \aand p\in\MFe \st\\[.1cm]
            \begin{array}{llllcllll}
                \aFe(\vulim+\vu^*,\vv)- \bFe(\vv,p) &=& \vg(\vv) &\quad \forall \vv\in\WFe\\[.1cm]
                 \bFe(\vulim,r)             &=& 0 &\quad \forall r\in\MFe.\\
            \end{array}
            \end{array}\]
Just as in Proposition \ref{prop:prop1}, this problem has a unique solution, as one can show, just as in Lemma  \ref{AFeContSMono}, that 
\[
\ab^*(\vu,\vv) := \ab(\vu+\vu^*,\vv)
\]
defines a continuous operator, strictly monotone on $\WFe$ and uniformly monotone on
$\tWFe$. 
Then, due to the bilinearity of $\bFe$,
$\vu := \vulim+\vu^*$, together with $p$ is a solution of    $\PFe $.

{{To show uniqueness}} we refer to Lemma \ref{unique} in Appendix \ref{ApUniqueness}.
\eproof
%
%
\section{Darcy as a limit of \DF - Simple Domain}\label{SimpDom}  
%
%
Suppose here that $\dom$ is a bounded domain in $\R^d$ with boundary $\bord$.
The object of this section is to show that the solution of the Darcy problem
$$\begin{array}{lllll}
	\alo\vu=-\grad p 	&\inn\dom\\
	\!\!\div \vu = \qo 	&\inn\dom\\
	\pp = \pbo  		&\on\bord
\end{array}$$
may be obtained as the limit of a sequence of solutions of the Darcy-Forchheimer problems
$$\begin{array}{lllll}
	\alo\vubet +\beto|\vubet|\vubet=- \grad \pbet 	&\inn\dom\\
	\!\!\div \vubet =\qo	&\inn\dom\\
	\pbet = \pbo 		&\on\bord,
\end{array}$$
as $\beto\rightarrow 0$.
As before we assume that  the tensor coefficient function
$\;\alo\,:\,\dom \longrightarrow R^{d,d},$ is such that
\begin{equation}\label{hypcoefsimp}\begin{array}{lllllll}
	\ulalo  |\vx|^2 &\leq& \vx\cdot\alo(\vy) \vx &\leq&\olalo |\vx|^2
		&\quad \forall \vy\in\dom,&\;\vx\in R^d,
\end{array}\end{equation}
and the coefficient $\beto$ of the nonlinear term is assumed to be a positive real parameter as we are merely interested in obtaining the Darcy problem as a limit of Forchheimer
problems.
Let
$$\begin{array}{lllll}
	\tY= H(\div,\dom)\qquad\qquad &\tX= L^2(\dom)\\[.2cm]
	\Y=W^3(\div,\dom)\qquad\qquad & \X= L^{\thalf}(\dom),
\end{array}$$
and recall that the image of the normal trace map on $\tY$ is $H^{-\half}(\bord)$
while the image of the normal trace map on $\Y$ is $W^{-\third,3}(\bord)$.
Also as before the data functions are assumed to be such that $\qq\in L^3(\dom)$
	and $\pbo\in W^{\third,\thalf}(\bord)\cap W^{\half,2}(\bord)$.

Define the bilinear forms $\ao$ and $\bo$ by
$$\begin{array}{rcccrcclc}
	\ao:\tY\times\tY &\longrightarrow &\R
		& \quad\aand \quad &
		 \bo:\tY\times\tX& \longrightarrow &\R  \\[.2cm]
	(\vu,\vv)\;\,& \mapsto&\!\! \dintx{\dom}{\alo \vu\cdot\vv}
		& & (\vv,\rr)\; \,&\mapsto&\!\! \!\!\dintx{\dom}{\!\!\div(\vv)\,\rr} ,
\end{array}$$
and the linear forms $\gg\in \tYp$ and $\ff \in \tXp$ by
$$\begin{array}{rcccrcclc}
	\gg\;:\:\tY &\longrightarrow &\hspace*{-2.5cm}\R & \quad\mbox{ and } \quad & \ff\;:\:\tX& \longrightarrow &\R
		\\[.2cm]
	\vv\;\,& \mapsto& \dpairbO& & \rr\; \,&\mapsto&\!\!\!\! \dintx{\dom}{\qo\rr},
\end{array}$$ 
so that the problem $\PDar$ can be written as
$$\PDar \quad\begin{array}{llll}
	\mbox{Find }\vu\in\tY \mbox{ and }\pp\in \tX \mbox{ such that }\\[.1cm]
	\begin{array}{llll}
		\ao(\vu,\vv) - \bo(\vv,\pp) &=& \gg(\vv)  &\quad\forall \vv\in \tY\\[.1cm]
		\bo(\vu,\rr) &=& \ff(\rr) 	&\quad\forall \rr\in \tX .
	\end{array}
\end{array}$$	
Since $\ao$ is elliptic (coercive) on the subset
${\tltY} =\{ \vv\in\tY\,:\,\bo(\vv,\rr)=0,\,\forall \rr\in\tX\}$
 and $\bo$ satisfies the inf-sup condition on $\tY\times\tX$:
 \[
 	\ulalo\|\vv\|_{\tY}\leq\ao(\vv,\vv)\;\forall\vv\in {\tltY} \qquad\aand\qquad
	\theto \|\rr\|_{\tX} \leq \dsup{\vv\in\tY}{\dfrac{\bo(\vv,\rr)}{\|\vv\|_{\tY}}},\;\forall \rr\in\tX,
\]
the Darcy problem $\PDar$ has a unique solution $(\vuo,\po)\in\tY\times\tX,$\cite{Brezzi[74]}.

To give the weak formulation of the Forchheimer problem note that since $\Y\subset\tY$ the bilinear form $\ao$ is also defined on $\Y\times\Y$ and that the
bilinear form $\bo$ is also defined on $\Y\times\X$ (even though $\X\not\subset\tX$). Further
$\bo$ also satisfies the analogous inf-sup condition on$\Y\times\X$ for some constant $\thetbo$;
see \cite{Summ-th[01]} or the more general version in  Lemma \ref{LeminfsupFe}.
Now define  the mapping $\abo$, linear in its second variable, by
$$\begin{array}{rcccrcccc}
		\abo\;:\:\Y\times\Y &\longrightarrow &\R \\[.2cm]
		(\vu,\vv)\;\,& \mapsto& \dintx{\dom}{(\alo+\beto |\vu|)\vu\cdot\vv} ,
\end{array}$$
and note that due to the regularity requirements on the data functions  $\pbo$ and  $\qo$ that
the linear forms $\gg$ and $\ff$ are defined and continuous on $\Y$ and $\X$, respectively,
(as well as on $\tY$ and $\tX$),
so that the problem $\PForch$ can be written as
$$\PForch \quad\begin{array}{llll}
	\mbox{Find }\vubet\in\Y \mbox{ and }\pbet\in \X \mbox{ such that }\\[.1cm]
	\begin{array}{llll}
		\abo(\vubet,\vv) - \bo(\vv,\pbet) &=& \gg(\vv)  &\quad\forall \vv\in \Y\\[.1cm]
		\bo(\vubet,\rr) &=&\ff(\rr) 	&\quad\forall \rr\in \X .
	\end{array}
\end{array}$$	
It is shown in \cite{Knabner-Summ} that the form $\abo$ is continuous, strictly monotone on $\Y$, and coercive
on ${\tlY} =\{ \vv\in\Y\,:\,\bo(\vv,\rr)=0,\,\forall \rr\in\X\}$ \cite[Proposition 1.2]{Knabner-Summ}
and that the Forchheimer problem $\PForch$ has a unique solution
$(\vubo,\pbeto)\in\Y\times\X,$ \cite[Theorem 1.8]{Knabner-Summ}.  Again see the more general vesion of this
reasoning in Lemma \ref{AFeContSMono}.
	
The demonstration that the solutions of the problems $\PForch$ converge to the solution
of $\PDar$ is based on a priori bounds for $\vubo\aand\pbeto$ independent of the
parameter $\beta$.   In this section we will drop the spaces in the notation for the norms
as only $\dom$ or $\bord$ appears.
\begin{lemma}\label{boundsimpdom}
There is a constant $C$ independent of $\beta$ such that for $\beta$ sufficiently small
\[
	\|\pbeto\|_{\X} + \|\vubo\|_{\tY}+\beta^{\third}\|\vubo\|_{0,3} \leq C.
\]
In addition,
\[
\beta\|\vubo\|_{0,3}\longrightarrow 0,\quad\mbox{ as }\quad \beta\longrightarrow 0.
\]
\end{lemma}
\bproof Taking $\vubo$ for
the test function $\vv$ in the first equation of $\PForch$ and noting that
$\Y\subset\tY,$  as in Section \ref{ExUniqDFFrac+DFMatrix} (cf. estimate (\ref{monotoneAFe}))
one obtains
$$\begin{array}{lllll}
	\ulalo\| \vubo\|^2_{0,2} + \tbeta \|\vubo\|^3_{0,3}
		&\leq&  \gg(\vubo) + b(\vubo,\pbeto)\\[.2cm]
		&\leq&  \|\gg  \|_{\tYp}  \|\vubo\|_{\tY}
			+\| \div(\vubo)\|_{0,3} \|\pbeto\|_{0,\thalf}\\[.2cm]
		&\leq&  \gORpbord \left(    \|\vubo\|_{0,2} + \|\div(\vubo)\|_{0,2}   \right)\\[.2cm]
			&&+ \; \| \div(\vubo)\|_{0,3} \|\pbeto\|_{0,\thalf}.
\end{array}$$
Next directly from the second equation of  $\PForch$
(regarded as an equation in $\MFe(\dom)^\prime=L^3(\dom)$),
we obtain
\begin{equation}\label{estLthreedivu}\begin{array}{lllll}
	\| \div(\vubo)\|_{0,3} =  \|\ff\|_{0,3}  ,
\end{array}\end{equation}
and, as there is a continuous embedding $\tX\emb \X$, i. e. $L^2(\dom)\emb L^{\thalf}(\dom)$
so that the second equation of
$\PForch$ holds for test functions in $L^2(\dom)$ as well as for those in $L^{\thalf}(\dom),$
we also have
\begin{equation}\label{estLtwodivu}\begin{array}{lllll}
	\| \div(\vubo)\|_{0,2} = \|\ff\|_{0,2}   .
\end{array}\end{equation}

Combining these last three inequalities we obtain
$$\begin{array}{lllll}
	\ulalo\| \vubo\|^2_{0,2} + \dfrac{c_\ell^3}{2} \beta \|\vubo\|^3_{0,3}
		&\leq& \gORpbord(\|\vubo\|_{0,2}+\|\ff\|_{0,2})
		   + \|\ff\|_{0,3}  \|\pbeto\|_{0,\thalf}
\end{array}$$
and
\begin{equation}\begin{array}{lllll}\label{estubet}
	{\half}\ulalo\| \vubo\|^2_{0,2} +\dfrac{c_\ell^3}{2} \beta \|\vubo\|^3_{0,3}
		&\leq& {\dfrac{1}{2\ulalo}} \gORpbord^2
		      +     \gORpbord\|f\|_{0,2}
		       + \|\ff\|_{0,3} \|\pbeto\|_{0,\thalf}\\[.3cm]
		 &\leq& D_1+\norm{0,3}{\ff}\norm{0,\thalf}{\pbeto},
\end{array}\end{equation}
where $D_1$ is a constant depending only on the coefficient $\alo$ and the data functions determining $\gg$ and $\ff$.
Then with the first equation of $\PForch$, we obtain, $\forall \vv\in\Y,$
$$\begin{array}{lllll}
	 |\bo(\vv,\pbeto)|&\leq& |\abo(\vubo,\vv)|+|\gg(\vv)|\\[.2cm]
	 	&\leq& \olalo\|\vubo\|_{0,2}\|\vv\|_{0,2}  +  \hbeta\normp{0,3}{\vubo}{2}\norm{0,3}{\vv}
				+\gORpbord\norm{\Y}{\vv}\\[.2cm]
		&\leq& \left(\olalo C_L\norm{0,2}{\vubo}  +  \hbeta\normp{0,3}{\vubo}{2}
				+ \norm{\third,\thalf}{\pbo}\right)\norm{\Y}{\vv},
\end{array}$$
where we again use $C_L$, respectively $C_\ell$, here specifically for the continuity constant for the embedding $L^3(\dom)\emb L^2(\dom)$, respectively $\ell^3(\R^d)\emb \ell^2(\R^d)$.
Using the inf-sup condition for $b$ on $\Y\times\X$ we have
$$
	\thetbo \|\rr\|_{0,\thalf} \leq \dsup{\vv\in \Y} {\dfrac{\bo(\vv,\rr)}{\|\vv\|_{\Y}}}.
$$
and thus
\begin{equation}\label{estpLthreehalves}
	\thetbo \|\pbeto\|_{0,\thalf}\;\leq\; \left(\olalo C_L\norm{0,2}{\vubo}  +  \hbeta\normp{0,3}{\vubo}{2} + \norm{\third,\thalf}{\pbo}\right).
\end{equation}
Plugging this estimate for $\pbet$ into (\ref{estubet}) we obtain
$$\begin{array}{lllll}
	\half\ulalo\| \vubo\|^2_{0,2} + \tbeta \|\vubo\|^3_{0,3}\\
		\qquad\qquad\leq D_1 + \dfrac{1}{\thetbo}\norm{0,3}{\ff}\left(\olalo C_L\norm{0,2}{\vubo}
			+  \hbeta\normp{0,3}{\vubo}{2} + \norm{\third,\thalf}{\pbo}\right).
\end{array}$$
Now using the inequality
$$\begin{array}{lllll}
	\dfrac{\olalo C_L}{\thetbo}\norm{0,3}{\ff}\norm{0,2}{\vubo}
		& \leq & \dfrac{4}{\ulalo}\left(\dfrac{\olalo C_L}{\thetbo}\right)^2\normp{0,3}{\ff}{2}
			+  \dps\fourth\ulalo\normp{0,2}{\vubo}{2} \\[.3cm]
		& \leq & D_2 +  \dps\fourth\ulalo\normp{0,2}{\vubo}{2}
\end{array}$$
it is easy to see that
$$\begin{array}{lllll}
	\dps\fourth\ulalo\| \vubo\|^2_{0,2} + \tbeta \|\vubo\|^3_{0,3}
		&\leq& D_1 + \dfrac{4}{\ulalo}\left(\dfrac{\olalo C_L}{\thetbo}\right)^2\normp{0,3}{\ff}{2}
			+\dfrac{1}{\thetbo}\norm{0,3}{\ff}\hbeta\normp{0,3}{\vubo}{2} \\[.4cm]&&\qquad
			+ \dfrac{1}{\thetbo}\norm{0,3}{\ff}\norm{\third,\tthird}{\pbo}\\[.4cm]
		&\leq& D_1 + D_2  +   \Cfour \beta \normp{0,3}{\vubo}{2} +D_3,
\end{array}$$
with constant terms $D_2$, which depends on $\ff, C_L,\olalo,\ulalo$ and $\thetbo$, and $D_3$, which
depends on $\ff,\pbo$ and $\thetbo$, and a constant coefficient $\Cfour$, which depends on $\ff$,
$\thetbo$ and $C_\ell$.
Now using Young's inequality, (if $p>0$ and $\frac{1}{p}+\frac{1}{q}=1$ then
$ab\leq \frac{a^p}{p}+\frac{b^q}{q}$) with $a=(\beta^{s}\norm{0,3}{\vubo})^2$, $b=1$, $p=\thalf$
and $q=3$ one obtains
$$\begin{array}{lllll}
	\dps{\fourth\ulalo\| \vubo\|^2_{0,2} + \tbeta \|\vubo\|^3_{0,3}}
	&\leq&\dps{ D_1 + D_2 + D_3 +\tthird \Cfour \beta^{1-2s}(\beta^s\norm{0,3}{\vubo})^3
			+ \third \Cfour\beta^{1-2s}}
\end{array}$$	
and that
$$\begin{array}{lllll}
	\dps{\fourth\ulalo\| \vubo\|^2_{0,2}  + \left( \gbeta\beta^{1-3s}- \tthird \Cfour \beta^{1-2s}   \right)
	(\beta^s\norm{0,3}{\vubo})^3 }
	&\leq&\dps{ D_1 + D_2 + D_3 + \third \Cfour\beta^{1-2s}}
\end{array}$$	
or, in particular,  that (with $s=\half$)
\begin{equation}\label{ineqntoestubeta}\begin{array}{lllll}
	\dps{\fourth\ulalo\| \vubo\|^2_{0,2} + \left( \gbeta\beta^{-\half}- \tthird \Cfour    \right)
	(\beta^{\half}\norm{0,3}{\vubo})^3 }
	&\leq&\dps{ D_1 + D_2 + D_3 + \third \Cfour}\deff \Dfour.
\end{array}\end{equation}	

Thus for $\beta$ sufficiently small, we obtain an a priori bound on $\ulalo^{\half}\| \vubo\|_{0,2}$:
\begin{equation}\label{estuLtwo}
	\dps{\ulalo^{\half}\| \vubo\|_{0,2}} \leq 2
		\left( \Dfour \right)^{\half},
\end{equation}
and also that
$$  
	(\beta^{\half}\norm{0,3}{\vubo})^3  \leq   \left(\gbeta \beta^{-\half}- \tthird \Cfour    \right)^{-1}
		\Dfour
$$  
so that
\begin{equation}\label{estBetahalfULthree}
	\beta^{\half}\norm{0,3}{\vubo} \longrightarrow 0 \qquad \mbox{ as }
		\qquad \beta  \longrightarrow 0.
\end{equation}
Rewriting (\ref{ineqntoestubeta}) as
$$\begin{array}{lllll}
	\dps{\fourth\ulalo\| \vubo\|^2_{0,2} + \gbeta\beta \|\vubo\|^3_{0,3}}
	&\leq&\dps{ \Dfour+\tthird \Cfour (\beta^{\half}\norm{0,3}{\vubo})^3 },
\end{array}$$	
we obtain in turn an a priori bound for $\beta^{\third} \norm{0,3}{\vubo}$:
\begin{equation}\label{estuLthree}
	\beta^{\third} \norm{0,3}{\vubo}  \leq
		\left( \dps{ \Dfour +\epsilon} \right)^{\third},
\end{equation}
with $\epsilon>0$ arbitrarily small for $\beta\leq\bar{\beta}_\epsilon$
for some $\bar{\beta}_\epsilon>0.$
Now combining  (\ref{estpLthreehalves}), (\ref{estuLtwo}) and (\ref{estuLthree}) one obtains
the following a priori bound on $\pbeto$ in $L^{\thalf}(\dom)$:
\begin{equation}\label{estp}\begin{array}{lll}
	\norm{0,\thalf}{\pbeto}
	&\leq &\dfrac{2 \olalo C_L}{\thetbo \ulalo^{\half}}\left( \dps{ \Dfour} \right)^{\half}
	+     \dfrac{C_L^3}{\thetbo}\beta^{\third}\left( \dps{ \Dfour +\epsilon} \right)^{\tthird}
	+\dfrac{1}{\thetbo}\norm{\third,\thalf}{\pbo}.
\end{array}\end{equation}
With (\ref{estp}), (\ref{estuLtwo}) and (\ref{estLtwodivu}) the lemma is completed.
\eproof

From (\ref{estuLtwo})  and (\ref{estp}), we conclude that
if $\{\beta_j\}$ is a sequence converging to 0 then there is a subsequence still denoted $\{\beta_j\}$ such that the sequences $\{\vubjo\}$ and $\{\pbetjo\}$
are weakly convergent in $(L^2(\dom))^d$ and in $L^{\thalf}(\dom)$, respectively:
\begin{equation}\label{wklims}
		\vubjo \wklim \vulim \mbox{ in } (L^2(\dom))^d\qquad \mbox{ and }
		\qquad \pbetjo \wklim \plim \mbox{ in } L^{\thalf}(\dom) ,
\end{equation}		
i. e. explicitly
\begin{equation}\label{IntConv}
	\dintx{\dom}{\alpha \vubjo\cdot \vv} \rightarrow \dintx{\dom}{\alpha \vulim \cdot \vv}
	\;\forall \vv\in (L^2(\dom))^{n}
	 \quad \mbox{ and }\quad
	 \dintx{\dom}{\pbetjo  q} \rightarrow  \dintx{\dom}{\plim  q}\;\forall q\in L^3(\dom) .
 \end{equation}
Further, (\ref{estBetahalfULthree}) implies that $\{\beta_{j}^{\half}\vubjo\}$
converges strongly to 0 in $L^3(\dom).$  Thus
\begin{equation}\label{ConvbetahalfUtozero}\begin{array}{lll}
	\left| \dintx{\dom}{\beta_j|\vubjo|\vubjo\cdot\vv}\right|
	\leq  \norm{0,3}{\vv}\norm{0,\thalf}{\beta_j|\vubjo| \vubjo}
	=      \norm{0,3}{\vv}\left( \dintx{\dom}{\beta_{j}^{\thalf}|\vubjo|^{3} }\right)^{\tthird}
	\\[.3cm]  \hspace{3.6cm}
	=       \norm{0,3}{\vv}\normp{0,3}{\beta_{j}^{\half}\vubjo}{2}
	\rightarrow 0  \mbox{ as } \beta_j \rightarrow 0  .
\end{array}\end{equation}

\begin{lemma} Assume that the spatial dimension $d$ satisfies $d\leq 6$.
Then the pair $(\vulim,\plim)$ defined by (\ref{wklims}) is a solution to $\PDar$ and hence
is the unique solution of $\PDar$: $\vulim = \vuo$ and $\plim=\po$.
\end{lemma}

\bproof A priori, $\vulim\in\Y\subset\tY$ and $\plim\in\X\not{\!\!\subset}\;\tX.$
It follows from (\ref{ConvbetahalfUtozero}) and (\ref{IntConv}) that
\begin{equation}\label{equoneDarY}\begin{array}{llll}
	\dintx{\dom}{\alo \vulim\cdot\vv} - \dintx{\dom}{\div(\vv)\plim}
		 =- \dpairbO &\forall \vv\in \Y.\\
\end{array}\end{equation}					
However, if $\vv\in \D(\dom),$
then $\dintx{\dom}{\!\!\alo \vulim\cdot\vv} - \dintx{\dom}{\!\!\div(\vv)\plim} =0 $ so that
$\grad\plim=-\alo\vulim \in (L^{3}(\dom))^d,$ and thus $\plim\in W^{1,\thalf} (\dom)$.
 By the Sobolev embedding theorem, \cite{Adams[75]} , we have, if $d\leq 6,$ then  $W^{1,\thalf}(\dom)\subset L^2(\dom)=\tX.$		
 Also we have supposed that $\pbo$ belongs to
$W^{\half,2}(\bord)$ as well as to $W^{\third,\thalf}(\bord)$.
Thus each of the terms of (\ref{equoneDarY}) is well defined for $\vv\in\tY$, and we have since
$\Y$ is dense in $\tY$  that
\begin{equation}\label{equoneDartY}\begin{array}{llll}
	\dintx{\dom}{\alo \vulim\cdot\vv} - \dintx{\dom}{\div(\vv)\plim}
		 =- \dpairbO &\forall \vv\in \tY.
\end{array}\end{equation}

Turning now to the second equation of $\PDar$, we recall that $\ffo=\qo$  belongs to $L^3(\dom)=\X^{\prime}$,
and thus also to $L^2(\dom)=\tX^{\prime}$.  As we have seen, the second equation of $\PForch$
implies that for each $\beta>0$, $\div(\vubo)=f \in L^3(\dom)$.
This with (\ref{estuLtwo}) implies that $\vubo$ is bounded in the $\tY$ norm and that for a subsequence							
 $\{\beta_\ell\}$ of $\{\beta_j\}, \vu_{\beta_{\ell}} $ converges weakly to $\vulim$ in $\tY$.
 It follows that
 $$\begin{array}{lll}\dintx{\dom}{\div(\vulim)\rr} = \dint{\dom}{\qo\rr} &\forall \rr\in \tX.
	\end{array}$$
Thus the pair $(\vulim,\plim)$ in $\tY\times\tX$ is a solution of $\PDar$ and
$(\vulim,\plim)= (\vuo,\po)$ by uniqueness.
\eproof
%
%
\section{Darcy as a limit of \DF - Domain with a Fracture}\label{ExUniqDFFrac+DMatrix}  
%
%
The object of this section is to obtain the original problem $\Prob$ (with Darcy flow in the
subdomains $\Omega_1$ and $\Omega_2$ but Forchheimer flow in the fracture $\gamma$)
as the limit of the problem $\PFe$ (with Forcheimer flow in the subdomains and in the
fracture) studied in Section \ref{ExUniqDFFrac+DFMatrix}
when the Forchheimer coefficient in the subdomains $\beta$ decreases to 0. In this section, as in  Section \ref{SimpDom}, for simplicity we shall
assume  that $\beta_i$ is the same constant, positive, real parameter for $i=1$ and $i=2$:
$$\beta_1=\beta_2=\beta>0.$$
(The tensors $\betgam,\ali\aand\algam$ (\ref{hypcoef}), as well as $\kappa$ (\ref{hypkap}), remain as in Section \ref{ExUniqDFFrac+DFMatrix}.)
For each $\beta$ sufficiently small, let $(\vubet,\pbet)\in \WFe\times\MFe$
be the solution of $\PFe$.  We will derive a priori bounds on $(\vubet,\pbet)$ which are
independent of $\beta$,
thus obtaining a limit function which we shall show is a solution to $\Prob$.
\begin{lemma}\label{boundsfractureddom}
There is a constant $C$ independent of $\beta$, such that, for each $\beta$ sufficiently small,
	$$\norm{\vW}{\vubet}+\norm{\MFe}{\pbet}
			+\beta^{\third}\dsum{\norm{0,3,\Omega_i}{\vubeti}}\leq C.$$
	In addition, $$\beta^{\half}\dsum{\norm{0,3,\Omega_i}{\vubeti}}\lra 0, \mbox{ as } \beta\lra 0.$$
\end{lemma}
\bproof  The proof follows closely the lines of the proof of Lemma \ref{boundsimpdom}.
Taking for test function $\vv=\vubet$ in the first equation of $\PFe$, noting that $\vubet\in\vW$ and that
	$\vg\in\vWp$,  and letting $\Cxi$ denote $\ulkap\min(1,2\xi -1)$ we obtain
\[
\begin{array}{lll}  \dsum{\left(
	\ulali \normp{0,2}{\vubeti}{2} +\tbeta\normp{0,3}{\vubeti}{3}    \right) }
		+  \ulalg \normp{0,2}{\vubetg}{2} +\tbetgam\normp{0,3}{\vubetg}{3}
		+ \dsum{ \Cxi \normp{0,2}{\vubeti\cdot\vn}{2}}\\[.2cm]
	\qquad  \leq \aFe(\vubet,\vubet) = \vg(\vubet) +\bFe(\vubet,\pbet)\\[.3cm]
	\qquad  \leq \norm{\vW^\prime}{\vg} \norm{\vW}{\vubet}
		+\norm{\MFep}{\!\!\Div\vubet}\norm{\MFe}{\pbet}\\[.2cm]
	 \qquad  =   \norm{\vWp}{\vg}
	 	\biggl(\norm{\vV}{\vubet}+\norm{\Mp}{\!\!\Div\vubet}
	 		+	\dsum{\norm{0,2,\gamma}{\vubeti\cdot\vn}}
		\biggr) +\norm{\MFep}{\!\!\Div\vubet}\norm{\MFe}{\pbet},
\end{array}
\]
and from the second equation we have $\Div\vubet=f$ so that
\begin{equation}\label{LtwoLthreeEstDivu}
\begin{array}{lcl}
\norm{\MFe^\prime}{\!\!\Div\vubet}& = &\norm{\MFe^\prime}{f}\\[.3cm]
\norm{\M^\prime}{\!\!\Div\vubet}& = &\norm{\M^\prime}{f}.
\end{array}
\end{equation}
Combining these estimates, analogously to (\ref{estubet}) we obtain
\[
\begin{array}{lll}
\dsum{\left(
	\ulali \normp{0,2}{\vubeti}{2} +\tbeta\normp{0,3}{\vubeti}{3}    \right) }
	+  \ulalg \normp{0,2}{\vubetg}{2} +\tbetgam\normp{0,3}{\vubetg}{3}
	+ \dsum{ \Cxi \normp{0,2}{\vubeti\cdot\vn}{2}}\\[.2cm]
	\qquad\leq
	 	\norm{\vWp}{\vg}
	 	\biggl(\norm{\vV}{\vubet}+\norm{\Mp}{f}
	 		+	\dsum{\norm{0,2,\gamma}{\vubeti\cdot\vn}}
		\biggr) +\norm{\MFep}{f}\norm{\MFe}{\pbet}\\[.2cm]
	\qquad\leq
		\norm{\vWp}{\vg}
	 	\biggl(\dsum{\norm{0,2}{\vubeti}} +\norm{0,3,\gamma}{\vubetg}
	 		+	\dsum{\norm{0,2,\gamma}{\vubeti\cdot\vn}}
		\biggr) +\norm{\MFep}{f}\norm{\MFe}{\pbet}  +D_1,
\end{array}
\]
where the constant term $\Done$ depends on $\vg$ and on $f$.
Then, using Young's inequality we have
\[
\begin{array}{lll}
\dsum{\left(\half
	\ulali \normp{0,2}{\vubeti}{2} +\tbeta\normp{0,3}{\vubeti}{3}    \right) }
	+ {\ulalg \normp{0,2}{\vubetg}{2} +\tthird\tbetgam\normp{0,3}{\vubetg}{3}}
	+\half\dsum{ \Cxi \normp{0,2}{\vubeti\cdot\vn}{2}}\\[.3cm]
	\qquad\qquad\leq
	 \norm{\MFep}{f}\norm{\MFe}{\pbet}+\Done+\Dtwo,
\end{array}
\]
{where $\Dtwo$ depends on $\vg,\ulali,\ulbetg,\aand \Cxi$}.
The inf-sup condition for $\bFe : \vWFe\times\MFe \longrightarrow \R$
together with the first equation of $\PFe$ yields
\[\begin{array}{lll}
	  \thetab\norm{\MFe}{\pbet}
	&\leq&	\dsup{\vv\in\vWFe}{ \dfrac{b(\vv,\pbet)}{\norm{\vWFe}{\vv}}}\\[.2cm]
	&\leq&	\dsup{\vv\in\vWFe}{ \dfrac{\ab(\vubet,\vv)-\vg(\vv)}{\norm{\vWFe}{\vv}}},
\end{array}\]
and using (\ref{estNonLinTerm})
\[\begin{array}{lll}
|\ab(\vubet,\vv)-\vg(\vv)|   &\leq&  \dsum{\left({\olali\norm{0,2}{\vubeti}\norm{0,2}{\vvi}}
					+\hbeta\normp{0,3}{\vubeti}{2}\norm{0,3}{\vvi}\right)}
					\\[.2cm]
				&&\quad+\olkap\xi\biggl(\dsum{\norm{0,2,\gamma}{\vubeti\cdot\vn}  }\biggr)
						\biggl(\dsum{\norm{0,2,\gamma}{\vvi\cdot\vn}  }\biggr)
					\\[.4cm]
				&&\quad+ \olalg\norm{0,2,\gamma}{\vubetg}\norm{0,2,\gamma}{\vvg}
				+\Clgs \olbetg\normp{0,3,\gamma}{\vubetg}{2}\norm{0,3,\gamma}{\vvg}
				+{\|\vg\|\ssWFep}\norm{\vWFe}{\vv}.\\[.2cm]
\end{array}\]

Then combining the last two estimates, analogously to (\ref{estpLthreehalves}) we have
\begin{equation}\label{infsupest}  \begin{array}{llll}
  	 \thetab\norm{\MFe}{\pbet}
&\leq&     \dsum{\left({\Cali\norm{0,2}{\vubeti}}
					+\hbeta\normp{0,3}{\vubeti}{2}
					+\olkap\xi\norm{0,2,\gamma}{\vubeti\cdot\vn} \right)}\\[.3cm]
			&& \qquad+ \Calg\norm{0,2,\gamma}{\vubetg}
				+ \Clgs\olbetg\normp{0,3,\gamma}{\vubetg}{2}
				+  {\|\vg\|\ssWFep}  .
\end{array}\end{equation}
So
\[\begin{array}{lll}
\dsum{\left(\half
	\ulali \normp{0,2}{\vubeti}{2} +\tbeta\normp{0,3}{\vubeti}{3}    \right) }
	+ \ulalg \normp{0,2}{\vubetg}{2} +\dps{\frac{c^3_\ell}{3}}\ulbetg\normp{0,3}{\vubetg}{3}
	+\half\dsum{ \Cxi \normp{0,2}{\vubeti\cdot\vn}{2}}\\[.2cm]\qquad
	\leq
	\Done+\Dtwo+ \Dthree + \dfrac{\norm{\MFe^\prime}{f}}{\thetab}
			\biggl( \dsum{\left({\Cali\norm{0,2}{\vubeti}}
					+\hbeta\normp{0,3}{\vubeti}{2}
					+\olkap\xi\norm{0,2,\gamma}{\vubeti\cdot\vn} \right)}\\[.2cm]\hspace*{8cm}
			 +\Calg\norm{0,2,\gamma}{\vubetg}
				+\Clgs \olbetg\normp{0,3,\gamma}{\vubetg}{2}   \biggr),
\end{array}
\]
where $\Dthree$ depends on $\vg,\thetab,\aand f$.
Then using Young's inequality (three times with exponents 2 and 2 and twice with exponents
3 and $\thalf$) we obtain
\begin{equation}\label{apribdsu}\begin{array}{lll}
\dsum{\left(\fourth
	\ulali \normp{0,2}{\vubeti}{2}
	+ \tbeta\normp{0,3}{\vubeti}{3}     \right)}
	+ \half \ulalg \normp{0,2}{\vubetg}{2} +\dfrac{c_\ell^3}{3} \ulbetg
	\normp{0,3}{\vubetg}{3}
	+ {\dfrac{1}{4}} {\dsum{ \Cxi \normp{0,2}{\vubeti\cdot\vn}{2}}}   \\[.2cm]\qquad
	\leq
	\Done+\Dtwo+ \Dthree + \Dfour + \dfrac{2}{3}\Cfive \dsum{  (\beta^{\half}\norm{0,3}{\vubeti})^{3} },
\end{array}
\end{equation}
with $\Dfour$ depending on $\olali, \ulali, \olalg, \ulalg, \olbetg, \ulbetg, \olkap, \ulkap, \xi, \thetab, \aand f$,
and, analogously to (\ref{ineqntoestubeta}),
\[\begin{array}{lll}
\dsum{\left(\fourth
	\ulali \normp{0,2}{\vubeti}{2} +\left(\gbeta\beta^{-\half}-\dfrac{2}{3}\Cfive\right)
	\left(\beta^{\half}\norm{0,3}{\vubeti}    \right)^3
	+\fourth\Cxi \normp{0,2}{\vubeti\cdot\vn}{2}\right)} \\[.4cm]\qquad
	+ \fourth \ulalg \normp{0,2}{\vubetg}{2} +\dfrac{c_\ell^3}{3} \ulbetg
	\normp{0,3}{\vubetg}{3}   \\[.4cm]\qquad
	\leq
	D_1+D_2+ D_3 + D_4  .
\end{array}
\]
(Recall that $\thetab$ does not depend on $\beta.$) Hence
\[  \beta^{\half}\norm{0,3}{\vubeti} \longrightarrow 0\mbox{ as } \beta \longrightarrow 0, \quad i=1,2,\]
 as in (\ref{estBetahalfULthree}), and, as in (\ref{estuLtwo}), (\ref{estuLthree}) each of the terms
$  \norm{0,2}{\vubeti},\;
	  \norm{0,2,\gamma}{\vubeti\cdot\vn},\;
	  \norm{0,2,\gamma}{\vubetg},  \;
	 \norm{0,3,\gamma}{\vubetg} \;$ and $\;
	 \beta^{\third}\norm{0,3}{\vubeti}  $
is bounded by a positive constant $\Dfive$, depending on $\olali, \olalg, \olbetg, \olkap, \xi, f$ and
$\vg$ but independent of $\beta$:
\begin{equation}\label{multiUbound}
	\norm{0,2}{\vubeti} \;+\;
	 \norm{0,2,\gamma}{\vubeti\cdot\vn} \;+\;
	 \norm{0,2,\gamma}{\vubetg}   \;+\;
	 \norm{0,3,\gamma}{\vubetg}     \;+\;
	 \beta^{\third}\norm{0,3}{\vubeti}   \;\leq\;  \Dfive.
\end{equation}
Combining (\ref{LtwoLthreeEstDivu}) and (\ref{multiUbound}) yields an a priori bound on $\vubeti$ in the
$H(\div,\Omega_i)$-norm.  Equation (\ref{LtwoLthreeEstDivu}) also gives an a priori bound on
$\divt\vubetg -\STbetj$ in the $L^{3}(\gamma)$-norm, which completes the a priori bound of
$\norm{\vW}{\vubet}$.

To bound $\pbet$ we recall (\ref{infsupest})
\[   \begin{array}{llll}
  \thetab \left(\dsum{\norm{0,\thalf}{\pbeti}} +  \norm{0,\thalf,\gamma}{\pbetg}  \right)
&\leq&     \dsum{\left({C_L\olali\norm{0,2}{\vubeti}}
					+C^3_\ell\beta\normp{0,3}{\vubeti}{2}
					+\xikap\norm{0,2,\gamma}{\vubeti\cdot\vn} \right)}\\[.4cm]
			&& + C_L\olalg\norm{0,2,\gamma}{\vubetg}
				+ \Clis\olbetg\normp{0,3,\gamma}{\vubetg}{2}
				+  {\|\vg\|\ssWFep}
\end{array}\]
and obtain for a positive constant $\Dsix$, depending on $\olali, \olalg, \olbetg, \olkap, \xi, f, \vg$ and
$\thetab$ but independent of $\beta$ :
\begin{equation}\label{boundpbetaMFe} \begin{array}{llll}
\dsum{\norm{0,\thalf}{\pbeti}} +  \norm{0,\thalf,\gamma}{\pbetg}
\leq \Dsix,			
\end{array}\end{equation}
which gives the a priori bound on $\norm{\MFe}{\pbet}$.
\eproof

\begin{theorem}\label{unisolP} Suppose $d\leq 6$.
	There exists a unique solution $(\vu,\pp)\in\vW\times\M$ of problem $\Prob$, and $(\vu,\pp)$
	is a weak limit of solutions $(\vubet,\pbet)\in \vW\times\M$ in the sense made precise below.
\end{theorem}
\bproof
The proof follows from the error bounds (\ref{multiUbound}),(\ref{boundpbetaMFe}) obtained in Lemma \ref{boundsfractureddom}.
As the spaces $H(\div,\Omega_i), \; $   $L^{\thalf}(\Omega_i),  \;$   $L^3(\gamma),\;$  $L^{\thalf}(\gamma)$ and $L^{2}(\gamma)$ are reflexive Banach spaces they are sequentially weakly compact. Thus from (\ref{estuLtwo})  and (\ref{estp}), we conclude that
if $\{\beta_{\ell}\}$ is a sequence converging to 0 then there is a subsequence $\{\beta_j\}$ such that the sequences $\{\vu_{\beta_j,i}\}$, $\{\pp_{\beta_j,i}\}$, $\{\vu_{\beta_j,\gamma}\}$,
$\{\pp_{\beta_j,\gamma}\}$, $\{\vu_{\beta_j,i}\cdot\vn\}$ and $\divt\vu_{\beta_j,\gamma} -\STbetj$
are weakly convergent in $H(\div,\Omega_i)$, $L^{\thalf}(\Omega_i)$, $L^3(\gamma)$, $L^{\thalf}(\gamma)$,
$L^{2}(\gamma)$ and in $L^3(\gamma)$ respectively:
\[\begin{array}{rlllrlll}
	\vu_{\beta_j,i} & \wklim & \vuilim	& \inn H(\div,\Omega_i)&
	\pp_{\beta_j,i} & \wklim & \pilim			& \inn L^{\thalf}(\Omega_i) \\[.2cm]
	\vu_{\beta_j,\gamma} & \wklim & \vuglim		& \inn        L^3(\gamma)&
	\pp_{\beta_j,\gamma} & \wklim & \pglim		& \inn	L^{\thalf}(\gamma) \\[.2cm]
	\divt\vu_{\beta_j,\gamma} -\STbetj & \wklim &\divugminust						& \inn        L^3(\gamma)&
	\vu_{\beta_j,i}\cdot\vn & \wklim & \vuinclim	& \inn 	L^{2}(\gamma)
\end{array}	\]
and
\[
\beta^{\half}\vubeti \longrightarrow 0 \inn L^3(\Omega_i).
\]
We remark that  since $\|\vubeti\cdot\vn\|_{ H^{-\half}(\partial\Omega_i)}
	\leq  C \|\vubeti\|_{ H(\mbox{\scriptsize{div}},\Omega_i)}$  is bounded independently of
$\beta$  that $\vu_{\beta_j,i}\cdot\vn$ converges  weakly  to $\vuilim\cdot\vn$
in $H^{-\half}(\partial\Omega_i)$.  
Then since $\vu_{\beta_j,i}\cdot\vn$ converges weakly  to $\vuinclim$ in $L^2(\gamma)$,
we have $\vuinclim = \vuilim\cdot\vn$

We also note that $\divt\vu_{\beta_j,\gamma}\in L^2(\gamma)$ so that
$\vu_{\beta_j,\gamma}\in H(\divt,\gamma)$.  Further, $\|\divt\vu_{\beta_j,\gamma}\|_{L^2(\gamma)}$
 and thus $\|\vu_{\beta_j,\gamma}\|_{H(\ssdiv,\gamma)}$ is
 bounded independently of $\beta$ so that $\vu_{\beta_j,\gamma}$ converges weakly to $\vuglim$
 in $H(\divt,\gamma)$.
 Following the same lines of reasoning we conclude that
 $$    \hat{\vu}_\gamma=\divt\tilde{\vu}_{\gamma}-(\tilde{\vu}_1\cdot\vn-\tilde{\vu}_2\cdot\vn) . 	$$
 Thus with $\vulim=(\tilde{\vu}_1,\tilde{\vu}_2,\tilde{\vu}_{\gamma})$
 we have $\vulim\in \vW$, and it is clear that
 \[  b(\vulim,r)=<\Div\vulim,r>_{\MpM }= f(r),\qquad \forall r\in \M,\]
 i. e. the second equation of $\Prob$ is satisfied by $\vulim$.


For each $\beta>0$, the first equation of $\PFe$ is
 \[\begin{array}{lllll}
 \ab(\vubet,\vv) - b(\vv,\pbet) = \vg(\vv), \qquad\forall\vv\in\vWFe,\\[.4cm]
\displaystyle{\sum_{i=1}^{2}\,\int_{\Omega_i}(\ali +\beta |\vubeti|) \vubeti\cdot\vv_i\,dx
		+\int_{\gamma}(\algam+\betgam|\vubetg|)\vubetg\cdot\vvg \,\dsig}  \\ \hspace{3cm}
		+\displaystyle{\sum_{i=1}^{2}\int_{\gamma}\frac{1}{\kgamn}
			(\xi\vubeti\cdot\vn+\xibar\vu_{\beta,i+1}\cdot\vn)\vv_i\cdot\vn\,\dsig}
			- b(\vv,\pbet) = \vg(\vv), \qquad\forall\vv\in\vWFe,
\end{array}\]
or
 \[\begin{array}{lllll}
\displaystyle{\sum_{i=1}^{2}\,\int_{\Omega_i}\beta |\vubeti| \vubeti\cdot\vv_i\,dx
		+\int_{\gamma}\betgam|\vubetg|\vubetg\cdot\vvg \,\dsig}
		=  -\displaystyle{\sum_{i=1}^{2}\,\int_{\Omega_i}\ali  \vubeti\cdot\vv_i\,dx
		-\int_{\gamma}\algam\vubetg\cdot\vvg \,\dsig}\\\hspace{3cm}
		-\displaystyle{\sum_{i=1}^{2}\int_{\gamma}\frac{1}{\kgamn}
			(\xi\vubeti\cdot\vn+\xibar\vu_{\beta,i+1}\cdot\vn)\vv_i\cdot\vn\,\dsig}
	 	+ b(\vv,\pbet) + \vg(\vv), \qquad\forall\vv\in\vWFe.
\end{array}\]
Then taking the limit as $\beta$ goes to 0 we have, due to (\ref{estNonLinTerm}) and Lemma \ref{boundsfractureddom}	
 \[\begin{array}{lllll}
\displaystyle{\lim_{\beta\rightarrow 0}\int_{\gamma}\betgam|\vubetg|\vubetg\cdot\vvg \,\dsig}
	=  -\displaystyle{\sum_{i=1}^{2}\,\int_{\Omega_i}\ali  \vuilim\cdot\vv_i\,dx
		-\int_{\gamma}\algam\vuglim\cdot\vvg \,\dsig}\\\hspace{3cm}
		-\displaystyle{\sum_{i=1}^{2}\int_{\gamma}\frac{1}{\kgamn}
			(\xi\vuilim\cdot\vn+\xibar\vulim_{i+1}\cdot\vn)\vv_i\cdot\vn\,\dsig}
	 	+ b(\vv,\plim) + \vg(\vv), \qquad\forall\vv\in\WFe,

\end{array}\]
and in particular, for test functions $\vv\in({\D}({\Omega}_1))^d\times({\D}({\Omega}_2))^d
					\times\{0\},$
\[
\displaystyle{\sum_{i=1}^{2}\,\int_{\Omega_i}\ali  \vuilim\cdot\vv_i\,dx }
	-\displaystyle{\sum_{i=1}^{2}\,\int_{\Omega_i}\div  \vvi\,\pilim\,dx} =0.
\]
Thus $\grad\pilim=-\ali  \vuilim \in L^2(\Omega_i)$ and therefore $\pilim\in W_{1,\thalf}(\Omega_i)$.
From the Sobolev embedding theorem we then have $\pilim\in L^2(\Omega_i)$ (for $d\leq 6$)
which means that $\plim\in\M$.  Now from the density of $\vWFe$ in $\vW$ we conclude that
 \[\begin{array}{lllll}
\displaystyle{\lim_{\beta\rightarrow 0}\int_{\gamma}\betgam|\vubetg|\vubetg\cdot\vvg \,\dsig}
	=  -\displaystyle{\sum_{i=1}^{2}\,\int_{\Omega_i}\ali  \vuilim\cdot\vv_i\,dx
		-\int_{\gamma}\algam\vuglim\cdot\vvg \,\dsig}\\\hspace{3cm}
		-\displaystyle{\sum_{i=1}^{2}\int_{\gamma}\frac{1}{\kgamn}
			(\xi\vuilim\cdot\vn+\xibar\vulim_{i+1}\cdot\vn)\vv_i\cdot\vn\,\dsig}
	 	+ b(\vv,\plim) + \vg(\vv), \qquad\forall\vv\in\vW.

\end{array}\]
Now there remains to see that
\[
	\displaystyle{\lim_{\beta\rightarrow 0}\int_{\gamma}\betgam|\vubetg|\vubetg\cdot\vvg \,\dsig}
	=\displaystyle{\int_{\gamma}\betgam|\vuglim|\vuglim\cdot\vvg \,\dsig},
			\qquad\forall\vvg\in L^3(\gamma).
\]
Toward this end we define a mapping on $L^3(\gamma)\times L^3(\gamma)$
by
\[
	(\vwg,\vvg)\mapsto \displaystyle{\int_{\gamma}\betgam|\vwg|\vwg\cdot\vvg \,\dsig}
\]
and the associated  mapping $\Cmap :L^3(\gamma)\longrightarrow L^{\thalf}(\gamma) $.
That the mapping $\Cmap $ is monotone and continuous can be shown as in the
proof of Lemma \ref{AFeContSMono} where the monotonicity and continuity of $\AFe $ are shown.
Therefore $\Cmap$ maps weakly convergent sequences to convergent sequences; see \cite{Zeidler[90]}.
Thus, since $\vubetg\wklim\vuglim$,
we have that $\Cmap(\vubetg)\rightarrow \Cmap(\vuglim)$ in $L^3(\gamma)$ which now yields that
\[
	a(\vulim,\vv)-b(\vv,\plim)=\vg(\vv),\qquad \forall\vv\in\vW .
\]
Thus $(\vulim,\plim)\in\vW\times\M$ is a solution of $\Prob$.

As in (\ref{monotoneAFe}) we see that $A$ is strictly monotone on $\vW$.  Thus we can refer to Lemma \ref{unique} for uniqueness.
\eproof
%
%
\appendix
\section{Appendix}
%
%
In this appendix for the sake of completeness we give the definition and some basic properties of the
spaces $W^p(\div,\Omega),$ and we include the demonstrations of some lemmas needed in the previous
sections.
\subsection{The spaces $W^p(\div,\dom)$}\label{Wpdivspaces}
We recall the definition given in \cite{Knabner-Summ}, \cite{Summ-th[01]} of the spaces $W^p(\div,\dom)$ for $\dom\subset\R^d$ a bounded domain in $\R^d$ and $p\in\R$ a number with $1\leq p$:
\begin{equation}\label{wpdiv}
   W^{p}(\div,\dom)\deff\{\vv\in (L^p(\dom))^d: \div\vv\in L^p(\dom)\}
\end{equation}
with norm
\[
	\norm{W^{p}(\ssdiv,\dom)}{\vv} \deff \norm{L^p(\dom)}{\vv}  +  \norm{L^p(\dom)}{\div\vv}.
\]
As pointed out in \cite{Knabner-Summ} and in \cite{Summ-th[01]} it suffices to note that $W^p(\div,\dom)$ is a closed subset of
$(L^p(\dom))^d$ to see that $W^p(\div,\dom)$ is a reflexive Banach space.
Further, normal traces of elements of $W^p(\div,\dom)$ belong to $W^{-\frac{1}{p},p}(\bord)$.
\subsection{A general uniqueness result}  \label{ApUniqueness}
The object here is to show a uniqueness  result that is used in the proofs of Theorems \ref{unisolPFe} and \ref{unisolP}: that if $A$, $B$ and $B^t$ are the operators
associated with a mixed formulation and $A$ is strictly monotone and $B$ is surjective, then the mixed problem has no more than
one solution;  more precisely
\begin{lemma}\label{unique}
Let $X$ and $Y$ be Hilbert spaces and let $a$ be a form, linear in its second variable, on $Y\times Y$
and b a bilinear form on $Y\times X$ and let $A : Y\lra Y^\prime$
and $B : Y\lra X^\prime$ be the associated linear operators defined by
$<A(v),w>_{Y^\prime,Y} = a(v,w),\; \forall w\in Y$ and $<B(v),r>_{X^\prime,X} = b(v,r),\; \forall r\in X$, respectively.
 Suppose further that $f\in X^\prime$ and $ g \in Y^\prime.$  Then if $A$ is strictly monotone and $B$ is surjective,
 the problem
\begin{equation}\label{}
({\mathcal P})\;
	\begin{array}{l}
	\mbox{ Find }u\in Y \aand p\in X \st\\
	\begin{array}{lllll}
		a(u,v) -b(v,p)&=& g(v)&\forall   v\in Y\\
		b(u,r)    &=& f(r) &\forall r\in X
	\end{array}
\end{array}
\end{equation}
has at most one solution.
\end{lemma}
\bproof  Suppose that $(u,p) \aand (w,s)\in Y\times X$ are solutions to $({\mathcal P})$.  Then
$$\begin{array}{lllll}
	a(u,v)-a(w,v) - b(v, p-s) &=&0, &\forall v\in Y\\
	b(u-w, r)		&=&0, &\forall r\in X,
\end{array}$$
and taking as test functions $v=u-w \aand r=p-s$ we obtain $a(u,u-w) - a(w,u-w)=0$ or
$<A(u)-A(w),u-w>_{Y^\prime,Y} =0.$  Then, as $A$ is strictly monotone we have $u=w$.
To see that $p=s$ we suppose the contrary and use the surjectivity of $B$ to obtain an element $v\in Y$ with
${<B(v),p-s>_{X^\prime,X}  \neq 0}$.
However, as $u=w$ we have $0=a(u,v)-a(w,v)=b(v,p-s) = {<~B(v),p-s>_{X^\prime,X} }$ contradicting the choice of $v$.
\eproof

\subsection{Some inf-sup conditions}  \label{ApInfSup}

In this paragraph we give proofs of the fact that some of the bilinear operators considered in the text satisfy the inf-sup condition.

\normalfont
  \medskip
  \noindent{\bf Proof of Lemma \ref{LeminfsupFe}:} \hspace*{3pt}\ignorespaces
This proof is just as that for the problem with Darcy flow in the fracture and in the subdomain (see the proof of \cite[Theorem 4.1]{MartinJR[05]}),
only modified  for  the inf sup condition in the $W^3(\div)\times L^{\thalf}$ setting as in the proof of \cite[Lemma A.3]{Knabner-Summ}.
It clearly suffices to show that the induced mapping $\BFe = \Div :\WFe\longrightarrow \MFep$ is surjective and has a continuous
right inverse. Given an element $\psi=(\psi_1,\psi_2,\psi_\gamma)\in \MFep$, to construct an element
$\vv_\psi=(\vv_1,\vv_2,\vvg)\in\WFe$ with $\BFe \vv_\psi = \psi$ and
$\norm{\WFe}{\vv_\psi}\leq C\norm{\MFep}{\psi}$ one solves the auxiliary problem $\Delta \phi = \hat{\psi}\inn \Omega,\; \phi=0 \on\Gamma$,
with right hand side $\hat{\psi}\in L^{3}(\Omega)$, the function that agrees with $\psi_i$ on $\Omega_i$.
The solution $\phi$ is in $W^{2,3}(\Omega)$ if $\Omega$ is sufficiently regular (otherwise we solve the
same homogenous Dirichlet problem on a larger more regular domain and take the restriction of the solution to $\Omega$),
and $ \|\phi\|_{2,3,\Omega}\leq C \|\hat\psi\|_{0,3,\Omega}$ with a constant $C$ that depends only on $\Omega$.
Then  $\hat{\vv} \deff \nabla\phi\in (W^{1,3}(\Omega))^d \subset W^3(\div,\Omega)$ and $\div\hat{\vv}=\hat\psi\in L^3(\Omega)$
so that $\|\hat{\vv}\|_{W^3(\ssdiv,\Omega)} \leq (1+dC)\|\hat\psi\|_{0,3,\Omega}$.
We also have
$\vvi \deff \hat{\vv}_{|\Omi} \in (W^{1,3}(\Omi))^d \subset W^3(\div,\Omi)$ with $\div{\vvi}=\psi_i\in L^3(\Omi)$
so that
$$\|{\vvi}\|_{W^3(\ssdiv,\Omi)} \leq (1+dC)\|\hat\psi\|_{0,3,\Omega}
	\leq(1+dC)\dfrac{1}{c_\ell}(\|\psi_1\|_{0,3,\Omega_1}+\|\psi_2\|_{0,3,\Omega_2})
	\leq \widetilde{C}\norm{\MFep}{\psi}.$$
As $\vvi\in (W^{1,3}(\Omi))^d$, we have
$\vvi\cdot\vn_i \in W^{\tthird,3}(\partial\Omi)$, where $\vn_i$ is the exterieur unit normal vector on $\Omi$, and it follows that
$\vvi\cdot\vn_i \in L^3(\gamma)\subset L^2(\gamma)$ and
$$\norm{0,2,\gamma}{\vvi\cdot\vni}
	\leq \norm{\tthird,3,\partial\Omega_i}{\vvi\cdot\vni}
	\leq\|{\vvi}\|_{W^3(\ssdiv,\Omi)}
	\leq \widetilde{C}\norm{\MFep}{\psi}.$$
Thus the pair $(\vv_1,\vv_2)$ is suitable for the first two components of $\vv$.

To obtain the third component $\vvg$, note that as $\hat{\vv}\in  W^3(\div,\Omega)$ we have $\vv_1\cdot\vn_1+\vv_2\cdot\vn_2=
0$ on $\gamma$,
and thus the problem in the fracture domain $\gamma$ is decoupled from that in the subdomains $\Omega_1$ and $\Omega_2$.
One has only to define $\vvg \deff \gradt \phi_\gamma$ where $\phi_\gamma$ is the solution of
$\Delta \phi_\gamma = {\psi_\gamma}\inn \gamma,\; \phi_\gamma=0 \on\partial\gamma$.  It is straightforward to verify now that
$\vv_\psi=(\vv_1,\vv_2,\vv_\gamma)$ is a suitable antecedent for $\psi$ and that the mapping $\psi\mapsto\vv_\psi$
is continuous from $\MFep$ into $\WFe$.
\eproof
\begin{lemma}\label{leminfsupWbM}
The inf-sup condition holds for the bilinear form $b:\vW\times\M\lra\R$; i. e.
there exists $\theta\in\R$ such that for each $r\in\M$
\begin{equation}\label{inf-supFracDomD_DF}
	\displaystyle{\sup_{\vv\in\vW} \frac{b(\vv,r)}{\norm{\vW}{\vv}}\geq  \theta\norm{\M}{r}}.
\end{equation}
\end{lemma}
\bproof The proof of this lemma is just as that of Lemma \ref{LeminfsupFe} only the auxiliary problem in the subdomains is in the
$H(\div)\times L^2$ setting.  As the auxiliary problems in the subdomains and in the fracture domain decouple no difficulty
arises from the fact that one of these involves $H(\div)\times L^2$ while the other involves $W^3(\div)\times L^{\thalf}$.
\eproof

\noindent
\underline{Acknowledgement:} We are grateful to the anonymous referee, who pointed out a considerable simplification of the original proof of Theorem 1.
\bibliography{ForchFrac}
\bibliographystyle{abbrv}


\end{document}